\documentclass[12pt]{amsart}

\usepackage{environ, xcolor}
\NewEnviron{commentA}{}

\newcommand\Aoff{\RenewEnviron{commentA}{}}
\Aoff{} 

\usepackage{amsmath, amssymb}
\usepackage{array}
\usepackage[frame,cmtip,arrow,matrix,line,graph,curve]{xy}
\usepackage{graphpap, color, paralist, pstricks}
\usepackage[mathscr]{eucal}
\usepackage[pdftex]{graphicx}
\usepackage[pdftex,colorlinks,backref=page,citecolor=blue]{hyperref}
\usepackage{cleveref}
\usepackage{tikz-cd, verbatim, transparent}

\newtheorem{theorem}{Theorem}[section]
\newtheorem{proposition}[theorem]{Proposition}

\newtheorem{lemma}[theorem]{Lemma}

\theoremstyle{definition}
\newtheorem{definition}[theorem]{Definition}
\newtheorem{example}[theorem]{Example}
\newtheorem{claim}[theorem]{Claim}

\newtheorem{remark}[theorem]{Remark}

\usepackage[margin=1in]{geometry} 
\usepackage{amsmath,amssymb,mathtools, mathrsfs,esint,tikz-cd,verbatim}

\newcommand{\Z}{{\mathbb Z}}
\newcommand{\Q}{{\mathbb Q}}

\newcommand{\C}{{\mathbb C}}

\newcommand{\on}[1]{\operatorname{#1}}

\newcommand{\im}{{\on{im}}}

\newcommand{\Gal}{{\on{Gal}}}

\newcommand{\Spec}{{\on{Spec}}}

\newcommand{\Pic}{{\on{Pic}}}


\title{Weighted surfaces with maximal Picard number}
\author{Louis Esser}
\author{Jennifer Li}

\address{Department of Mathematics, Princeton University, Fine Hall, Washington Road, Princeton, NJ 08544-1000, USA}
\email{esserl@math.princeton.edu}

\email{jenniferli@princeton.edu}

\begin{document}

\begin{abstract}
    An algorithm due to Shioda computes the Picard number
    for certain surfaces which are defined by a single equation
    with exactly four monomials, called Delsarte surfaces. We 
    consider this method
    for surfaces in weighted projective $3$-space with quotient
    singularities.  We give a criterion for such a 
    weighted Delsarte surface $X$ to
    have maximal Picard number. This condition is surprisingly related to 
    the automorphism group of $X$. For every positive integer $s$, we find 
    a weighted Delsarte surface with geometric genus $s$ and maximal 
    Picard number. We show that these examples are elliptic surfaces, proving
    that elliptic surfaces of maximal Picard number and arbitrary geometric
    genus may be embedded as quasismooth hypersurfaces in weighted 
    projective space.
\end{abstract}

\maketitle

\section{Introduction}

The Picard number $\rho(X)$ of a smooth complex 
projective variety $X$ is the rank of the N\'{e}ron-Severi 
group $\mathrm{NS}(X)$ of algebraic cycles of codimension $1$
in $X$, up to algebraic equivalence.  In families of algebraic varieties,
it's often true that ``most" fibers have small $\rho$ and carry only
the algebraic classes required to exist via, e.g., their projective embedding.
On the other hand, special members of the family may have ``extra" 
algebraic classes.

The size that $\mathrm{NS}(X)$ can attain is bounded by cohomology;
in particular, we always have 
that $\rho$ is at most the Hodge number $h^{1,1}(X)$.  If $\rho = h^{1,1}$,
we say that $X$ has \textit{maximal Picard number}.

The computation of the Picard number for a specific $X$ 
is a very difficult problem in general, even when $X$ is a surface.  As a 
result, relatively few surfaces with maximal Picard number
are known to exist
(see \cite{Beauville} for a survey).  For instance, it remains an 
open question whether there exists a smooth surface in $\mathbb{P}^3$ with
degree $d \geq 7$ and maximal Picard number.

One example for which $\rho$ is explicitly computable is the Fermat
surface $Y_d \coloneqq \{y_0^d + y_1^d + y_2^d + y_3^d  = 0\} 
\subset \mathbb{P}^3$.
Leveraging the large group of symmetries of $Y_d$,
Shioda and Aoki were 
able to prove a formula for $\rho(Y_d)$ for every $d$ \cite{Shioda82,Aoki}.
Shioda later observed \cite{Shioda86} that one can similarly 
compute the Picard number when $X$ is a surface defined by a 
single polynomial
equation with exactly four monomials, satisfying certain conditions.
An $X$ of this type is often called a \textit{Delsarte surface},
after \cite{Delsarte}. 
More precisely, a Delsarte surface $X$ is
defined by an equation of the form
\begin{equation}
\label{eq:Delsarte_surf}
\sum_{i=0}^{3} x_{0}^{a_{i0}} x_{1}^{a_{i1}} x_{2}^{a_{i2}} x_{3}^{a_{i3}} = 0. 
\end{equation}
We may associate to $X$ the integral $4 \times 4$ matrix 
$A := (a_{ij})$. 
Shioda's method relies on expressing $X$
as a quotient of some Fermat surface, up to a birational transformation.

Subsequent work has utilized Shioda's algorithm to determine the Picard numbers
for all Delsarte surfaces in $\mathbb{P}^3$ with ADE singularities \cite{Heijne16}
and for Delsarte K3 surfaces in weighted projective space, with connections to 
Berglund-H\"{u}bsch-Krawitz (BHK) mirror symmetry \cite{Goto,Kelly,LO}.
In both of these special cases, all surfaces considered have canonical 
singularities.

In this paper, we consider the more general case of a
\textit{weighted Delsarte surface},
which is a Delsarte
surface in weighted
projective space $\mathbb{P} := \mathbb{P}(a_{0},a_1,a_2,a_3)$ with 
quotient singularities.  In this generality, many new phenomena arise.
For instance, it often happens that such a surface $X$ has 
ample canonical
class $K_X$, but a desingularization of $X$ 
has Kodaira dimension $\kappa < 2$; such an $X$ can even be rational.  
There are also many new examples with maximal Picard number, where
for a singular surface we take this to mean that a desingularization 
(equivalently, any desingularization) has maximal Picard number.

Our first main result is a criterion for a 
weighted Delsarte surface $X \subset \mathbb{P}$ to have 
maximal Picard number. This highlights a surprising connection with the toric
automorphism group of $X$, which is the subgroup of the torus action on $\mathbb{P}$
leaving $X$ 
invariant (see \Cref{def:toricAutomorphismGroup}). 
This group is easily computed from the exponent matrix $A$
from \eqref{eq:Delsarte_surf}.
In the following theorem, recall that
the exponent $e$ of a finite group $G$ is the least common multiple 
of orders of all elements in $G$.

\begin{theorem}
\label{thm:intro_aut}
Suppose that $X$ is a well-formed Delsarte surface in a weighted 
projective $3$-space $\mathbb{P}$ with at worst quotient singularities.
If the group 
$\mathrm{Aut}_{\mathrm{tor}}(X)$ 
of toric automorphisms of $X$ is a finite
group of exponent $e \in \{1,2,3,4,6\}$, then $X$ has maximal 
Picard number.  Further, if $e \leq 3$, then $X$ is rational.
\end{theorem}

This theorem can be viewed as an analog of the fact that the Fermat surface
$Y_d \subseteq \mathbb{P}^3$
has maximal Picard number if and only if $d = 1,2,3,4$, or $6$.  However,
there are many examples to which \Cref{thm:intro_aut} applies where $X$
is not a quotient of one of these small degree Fermat surfaces.
We note that since $\mathrm{Aut}_{\mathrm{tor}}(X) \subseteq \mathrm{Aut}(X)$,
it's of course sufficient for the conclusion of \Cref{thm:intro_aut} 
to assume that $\mathrm{Aut}(X)$ itself has exponent $1,2,3,4,$ or $6$.
This criterion applies to a diverse range of examples
(see Examples \ref{ex:rational_ex} through \ref{ex:gen_type_low_exp}).

The methods used to prove rationality in \Cref{thm:intro_aut}
also yield the following statement that does not require Delsarte assumptions.

\begin{theorem}
\label{thm:intro_rat}
Suppose $X$ is a well-formed quasismooth surface in weighted projective $3$-space.
If the geometric genus of $X$ equals $0$, then $X$ is rational.
\end{theorem}

This theorem follows from the toric geometric interpretations of geometric 
genus and Kodaira dimension, originally due, respectively, to Khovanskii 
\cite{Khovanskii} and Fine \cite{Fine} (see also \cite[Section 4, Appendix]
{Reid85}), as well as highly non-trivial classifications of three-dimensional
lattice polytopes with certain properties \cite{AKW,AWW,BKS,Treutlein}.  
Though \Cref{thm:intro_rat} is a relatively straightforward consequence of 
known results, we emphasize it here because it seems somewhat buried in the 
literature.  As recently as 2025, it was asked as an open question whether
there exists a non-rational quasismooth surface in weighted projective $3$-space 
satisfying 
$h^2(X,\mathcal{O}_X) = 0$ \cite[Question 3.20]{Chitayat}.  \Cref{thm:intro_rat}
answers this question in the negative (since $p_g(X) = h^2(X,\mathcal{O}_X)$ for 
$X$ quasismooth).

In ordinary projective space $\mathbb{P}^3$,
one of the main consequences 
of the classification in
\cite{Heijne16} is that there are only 
three Delsarte surfaces in $\mathbb{P}^3$ with quotient singularities,
up to isomorphism, 
with maximal Picard number and degree $d > 4$ (one with $d = 5$ and two 
with $d = 6$).
In particular, the geometric genus is bounded for examples of 
maximal Picard number: each example satisfies $p_g \leq 10$.
When $X$ is a surface with $p_g = 0$, the Lefschetz (1,1)-Theorem 
implies that $X$ has maximal Picard number (see \Cref{rem:pg0-maxPicardNumber}).
But for increasing $p_g$ it becomes harder to find examples.
In this paper, we show that, unlike in ordinary projective space,
surfaces with maximal Picard number and unbounded
$p_g$ exist in weighted projective $3$-space.
In fact, we find a quasismooth surface in weighted projective $3$-space
with maximal Picard number for every possible geometric genus.

\begin{theorem}
\label{thm:intro_elliptic}
For every positive integer $s$, there exists a well-formed
quasismooth Delsarte 
surface $X = X_s$
in weighted projective $3$-space with geometric genus $p_g = s$ 
and maximal Picard number.
\end{theorem}

One pleasant feature of these examples is the ease with which they
are defined.  For instance, we can define an $X_s$
for each $s$ satisfying the condition
of \Cref{thm:intro_elliptic}
by the remarkably simple equation:
$$x_0^{4s} + x_1^{4s} x_2 + x_2^{2s+1} x_3 + x_3^{2s+1} = 0.$$
We prove that this equation in fact defines an \textit{elliptic} surface,
so the desingularization of $X_s$ has Kodaira dimension $\kappa = 1$
when $s \geq 2$. We show that each of these elliptic
surfaces, despite
having maximal Picard rank, has Mordell-Weil rank $0$.

There are many
sequences of surfaces of this kind; in contrast, there appear to be only sporadic 
examples of Delsarte surfaces of maximal Picard number
with desingularization of general type.

Finally, as part of this project
we develop computational tools in \texttt{MAGMA}
for calculating the Lefschetz number $\lambda = b_2 - \rho$ of any 
weighted Delsarte
surface with quotient singularities
(though no results in the paper rely on this code). These tools are
publicly 
available online at \cite{magma}.

\Cref{sect:prelim} introduces Delsarte hypersurfaces in weighted
projective space and their properties, as well as the necessary
preliminaries on birational geometry and cohomology.  
\Cref{sect:algorithm} develops Shioda's algorithm in the weighted
setting; our description more heavily uses the 
language of toric geometry.  \Cref{sect:rat} introduces Newton polytopes
of hypersurfaces and their connections to birational invariants.  This
leads to the proof of \Cref{thm:intro_rat}.
In \Cref{sect:aut}, we prove 
\Cref{thm:intro_aut} and illustrate how this criterion
can be used to find different types of surfaces
with maximal Picard number.  Finally, 
in \Cref{sect:unbounded_genus}, we prove \Cref{thm:intro_elliptic}
and show that the surfaces we obtain are in fact elliptic.

\noindent \textit{Acknowledgements.} We thank 
Lena Ji, J\'{a}nos Koll\'{a}r, and Jos\'{e} Ignacio Y\'{a}\~{n}ez
for many helpful suggestions and comments.

\section{Preliminarites}
\label{sect:prelim}

\noindent \textbf{Notation:}
Throughout this paper, we work over the complex numbers $\C$.
For a positive integer $b$ and any integer $a$,
we use the notation $a \bmod b$ to refer to the
unique integer in $\{0,\ldots,b-1\}$ equivalent to
$a$ modulo $b$.

This section will introduce the necessary background
on Delsarte hypersurfaces and cohomology that will be 
used in subsequent sections.

\subsection{Introduction to Delsarte hypersurfaces}
\label{subsect:Delsartedef}

In this subsection, we'll see that various properties
of a Delsarte hypersurface may be calculated in terms of a 
square matrix $A$ associated to its equation.
These hypersurfaces were studied
in various other works, such as \cite{ABS}, where many of the
results of this subsection appear.  The study of Delsarte
hypersurfaces dates back to \cite{Delsarte}.

\begin{definition}
\label{def:delsartePolynomial}
A \textit{Delsarte polynomial} is a polynomial $f$ with 
the same number of monomials as variables.
\end{definition}

Given a Delsarte polynomial $f$ in $n+2$ variables $x_0,\ldots,x_{n+1}$,
we may write $f$ as
$$f = \sum_{i = 0}^{n+1} K_i \prod_{j = 0}^{n+1} x_j^{a_{ij}},$$
where we assume that all coefficients $K_i$ are nonzero.  The exponents in 
this equation form a square matrix $A = (a_{ij})$.  From now on, we will 
only consider the case
that $f$ is \textit{invertible},
that is, this matrix $A$ is invertible over the rational numbers $\Q$.

We may naturally view $f$ as the defining equation of a hypersurface in 
weighted projective space as follows.  For a fixed $j$, define $q_j$, the
$j$th \textit{charge} of $f$, to be the sum of the entries in the $j$th row 
of $A^{-1}$. Let $m$ be the least common denominator of $q_0,\ldots,q_{n+1}$, 
and  set $a_j \coloneqq m q_j$.  Then, 
if we assign the weights $a_j$ to the 
variables $x_j$, the equation $f$ is homogeneous of weighted degree $m$ 
(this can be seen from the matrix equation $AA^{-1} \mathbf{1} = \mathbf{1}$,
where $\mathbf{1}$ is the
vector of all $1$'s).  We further define the positive integer $d$ as the
least common denominator of all entries of $A^{-1}$, so that $m$
divides $d$. Later on, we'll see that this $d$ is the degree
of the Fermat hypersurface of which the Delsarte hypersurface
defined by $f$ is a (birational) quotient.

Given that $A$ is invertible,
the weights $a_0,\ldots,a_{n+1}$ are uniquely defined by the conditions that:
(1) $f$ is weighted homogeneous with variables of
weights $a_0,\ldots,a_{n+1}$, and (2)
$\gcd(a_0,\ldots,a_{n+1}) = 1$.  Thus, $f$ determines a 
hypersurface in $\mathbb{P}(a_0,\ldots,a_{n+1})$, the 
\textit{weighted projective space} with weights $a_0,\ldots,a_{n+1}$.
For more background on weighted projective spaces and their 
subvarieties, see \cite{Iano-Fletcher}.

We say that $\mathbb{P}(a_0,\ldots,a_{n+1})$ is \textit{well-formed}
if $\gcd(a_0,\ldots,\widehat{a_j},\ldots,a_{n+1}) = 1$ for all 
$j = 0,\ldots,n+1$.  Though the gcd of all the weights is always $1$
in our construction, this stronger well-formedness condition may not always
be true for a given Delsarte polynomial $f$.
However, every weighted projective space $\mathbb{P}$ is isomorphic as a variety
to a well-formed one, $\mathbb{P}^{\, \prime}$, via some coordinate substitution,
which will replace $f$ with another, slightly modified, Delsarte polynomial.
Therefore, we will only consider Delsarte polynomials that give rise
to well-formed weighted projective spaces.  We further define that 
a subvariety $X \subseteq \mathbb{P}(a_0,\ldots,a_{n+1})$ is \textit{well-formed}
if $\mathbb{P}$ is well-formed and 
$X \cap \mathbb{P}_{\mathrm{sing}}$ has 
codimension at least $2$ in $X$.
For this paper, we will always assume this property holds, as in the following
definition:

\begin{definition}
A \textit{Delsarte hypersurface} is a well-formed hypersurface $X_m$ 
in a weighted projective 
space $\mathbb{P}(a_0,\ldots,a_{n+1})$ defined by an invertible
Delsarte polynomial $f$ that is homogeneous of weighted degree $m$, 
where the weights $a_0,\ldots,a_{n+1}$ and degree $m$
are determined from $f$ via the above procedure.
\end{definition}

When we write $X_m$ to represent a hypersurface of degree $m$,
we may simply write $X$ if the context is clear.
For a normal, well-formed hypersurface $X_m = \{f = 0\}$
in the weighted projective space
$\mathbb{P}(a_0,\ldots,a_{n+1})$, the adjunction formula
$K_X = \mathcal{O}_X(m - a_0 - \cdots - a_{n+1})$ holds. 
We also say that $X$ is \textit{quasismooth} if the affine cone
$\{f = 0\} \subseteq \mathbb{A}^{n+2}$ is smooth away from the origin.
Quasismooth hypersurfaces have only cyclic quotient singularities.

We note that the isomorphism class of a Delsarte hypersurface $X$ is
independent of the choice of coefficients $K_i$ in the 
defining equation $f$,
so long as they are nonzero.  
This is because any two equations with different
coefficients are isomorphic via multiplication 
by some element of the torus.
Therefore, the matrix $A$ determines the properties of the 
hypersurface. From now on, 
we will omit the coefficients $K_i$, setting them
equal to $1$ for simplicity.

Next, we'll prove some basic facts about toric automorphisms of Delsarte 
hypersurfaces.

\begin{definition}
\label{def:toricAutomorphismGroup}
Let $f$ be a Delsarte polynomial and $X$ the 
corresponding hypersurface. The
\textit{toric automorphism group} $\mathrm{Aut}_{\mathrm{tor}}(f)$ of $f$ is 
the group of diagonal transformations $x_j \mapsto c_j x_j$ which leave the
polynomial $f$ invariant.  Similarly, we define
$\mathrm{Aut}_{\mathrm{tor}}(X)$ as the
group of diagonal automorphisms of $\mathbb{P}$ that leave $X$ invariant.
Equivalently, this is the subgroup of the torus action on $\mathbb{P}$
that stabilizes $X$.
\end{definition}

Any diagonal automorphism $x_j \mapsto c_j x_j$ of $f$ must leave every 
monomial of $f$ invariant individually, so
$$\prod_{j = 0}^{n+1} (c_j x_j)^{a_{ij}} = \prod_{j=0}^{n+1} x_j^{a_{ij}}$$
for each $i = 0,\ldots,n+1$.  This implies $\prod_{j = 0}^{n+1} c_j^{a_{ij}} = 1$;
in particular, after taking the logarithm of both sides, we see that
$$\begin{pmatrix}
    \log |c_0| \\ 
    \vdots \\ 
    \log |c_{n+1}| 
\end{pmatrix} \in \ker(A).$$
The kernel of $A$ is trivial, hence each $c_j$ has magnitude $1$, so that 
$c_j = e^{2 \pi \theta_j \sqrt{-1}}$ for some $\theta_j$. Denote by $\theta$
the vector $(\theta_j)$.  Now the identities
$\prod_{j = 0}^{n+1} c_j^{a_{ij}} = 1$ for each $i$ imply
$A \theta \in \Z^{n+2}$.  In particular, $\theta \in A^{-1} \Z^{n+2}$.
Each $\theta_j$ is only defined up to 
addition by any integer so there is a natural identification
$$\mathrm{Aut}_{\mathrm{tor}}(f) \cong A^{-1} \Z^{n+2}/ \Z^{n+2} \subseteq 
(\Q/\Z)^{n+2}.$$
In particular, $\mathrm{Aut}_{\mathrm{tor}}(f)$
is a finite abelian group of order $\left|\det(A)\right|$.

\begin{lemma}
\label{lem:Aut_tor}
There is an exact sequence
$$0 \rightarrow \Z/m \rightarrow \mathrm{Aut}_{\mathrm{tor}}(f) \rightarrow 
\mathrm{Aut}_{\mathrm{tor}}(X) \rightarrow 0,$$
where we view $\mathrm{Aut}_{\mathrm{tor}}(f)$ as a subgroup of $(\Q/\Z)^{n+2}$
as above and $\Z/m$ is the cyclic subgroup of 
order $m$ generated by the vector
of charges
$$
\begin{pmatrix}
    q_0 \\
    \vdots \\
    q_{n+1}
\end{pmatrix} = 
A^{-1} \begin{pmatrix}
    1 \\
    \vdots \\ 
    1
\end{pmatrix}.$$
\end{lemma}

\begin{proof}
Any element in $\mathrm{Aut}_{\mathrm{tor}}(X)$ lifts to a diagonal 
automorphism sending $f$ to a nonzero multiple.  After scaling appropriately,
we can arrange that this multiple is $1$, so the last map is surjective.

The toric automorphisms of $\mathbb{A}^{n+2}$ that descend to the trivial
automorphism on $\mathbb{P}$ are precisely those of the form
$(x_0,\ldots,x_{n+1}) \mapsto (t^{a_0} x_0,\ldots,t^{a_{n+1}} x_{n+1})$
for some $t \in \C^*$.
The elements of $(\Q/\Z)^{n+2}$ of this form are rational multiples of the 
vector of weights, which is in turn a multiple of the vector of charges.
The charge vector $A^{-1} \mathbf{1}$ is primitive in
$A^{-1} \Z^{n+2}$ so it generates the kernel of 
$\mathrm{Aut}_{\mathrm{tor}}(f) \rightarrow \mathrm{Aut}_{\mathrm{tor}}(X)$.
This subgroup has order $m$, the degree of the hypersurface, by our definition
of $m$.
\end{proof}

\begin{example}
\label{ex:Fermat}
The most fundamental example of a Delsarte polynomial is the \textit{Fermat
polynomial} $f = x_0^d + \cdots + x_{n+1}^d$ of degree $d$.  This has 
associated matrix $A = d I_{n+2}$, so, as expected, the weights are all
equal to $1$ and the degree is $m = d$.  Hence $f$ defines a hypersurface $X
\subseteq \mathbb{P}^{n+1}$ of degree $d$, called the
\textit{Fermat hypersurface}. The toric automorphism group of the Fermat
equation is $\mathrm{Aut}_{\mathrm{tor}}(f) = d^{-1} \Z^{n+2} / \Z^{n+2}
\cong (\Z/d)^{n+2}$, corresponding to transformations of each coordinate
by $d$th roots of unity, while 
$\mathrm{Aut}_{\mathrm{tor}}(X) = (\Z/d)^{n+2}/(\Z/d)$, the quotient of
this group by overall scalar multiplications.
\end{example}

\subsection{Delsarte hypersurfaces as rational quotients}
\label{subsect:quotient_setup}

One key property of a Delsarte hypersurface is that it is 
birational to a quotient
of a Fermat hypersurface in $\mathbb{P}^{n+1}$ by a finite toric
subgroup of its automorphism group.  This was stated by Shioda in 
\cite{Shioda86} for Delsarte surfaces in $\mathbb{P}^3$.
In this subsection, we give this quotient description in the more
general setting of weighted Delsarte surfaces using the language of toric
geometry; this description will be essential in later sections.

For this subsection, fix an irreducible Delsarte hypersurface 
$X \subseteq \mathbb{P} \coloneqq \mathbb{P}(a_0,\ldots,a_{n+1})$
with defining polynomial $f(x_0,\ldots,x_{n+1})$ and matrix $A$.
Consider the following sequence of maps of toric lattices $M', M,$ and $M''$
(each isomorphic to $\Z^{n+2}$):
    $$M' \xleftarrow{B^{\mathsf{T}}} M \xleftarrow{A^{\mathsf{T}}} M''.$$
Here $B = (b_{ij})$ is the integer matrix $d A^{-1}$.  Starting from these lattice maps,
we build the following diagram:

\begin{equation}
\label{eq:toricsetup}
\begin{tikzcd}
\Spec(\C[M']) \arrow[r] \arrow[d,hook] & 
\Spec(\C[M]) \arrow[r] \arrow[d,hook] &
\Spec(\C[M'']) \arrow[d,hook] \\
\mathbb{A}^{n+2}\setminus \{0\} \arrow[d] \arrow[r,dashed] & \mathbb{A}^{n+2}\setminus \{0\} \arrow[d] \arrow[r,dashed] & 
\mathbb{A}^{n+2}\setminus \{0\} \arrow[d] \\
\mathbb{P}^{n+1} \arrow[r,dashed] & \mathbb{P} \arrow[r,dashed] & \mathbb{P}^{n+1} \\
Y \arrow[r,dashed,"\varphi"] \arrow[u, hook] & X \arrow[r,dashed,"\psi"] \arrow[u, hook] & Z 
\arrow[u, hook].
\end{tikzcd}
\end{equation}

Here $Y$ and $Z$ are defined to be the Fermat hypersurfaces of degrees
$d$ and $1$, respectively.
The morphisms in the first row of \eqref{eq:toricsetup} are the 
homomorphisms of tori associated to the indicated lattice maps.  We use coordinates
$y_0,\ldots,y_{n+1}$ for $\C[M']$, $x_0,\ldots,x_{n+1}$ for $\C[M]$, and
$z_0,\ldots,z_{n+1}$ for $\C[M'']$.  In these
coordinates, the morphisms in
the first row correspond to ring homomorphisms 
$x_i \mapsto \prod_{j = 0}^{n+1} y_j^{b_{ij}}$ and 
$z_i \mapsto \prod_{j = 0}^{n+1} x_j^{a_{ij}}$, respectively.  
These induce rational maps on affine space (minus the origin) as in the 
second row.

We will next show that these rational 
maps descend to toric rational maps of the weighted projective spaces
indicated.  Indeed, consider the following commutative diagram
of cocharacter lattices for the toric varieties from the second and third
rows of
\eqref{eq:toricsetup}:
\begin{equation}
\label{eq:toricsetup_lat}
\begin{tikzcd}
N' \arrow[r,"B"] \arrow[d] & N \arrow[r,"A"] \arrow[d] & N'' \arrow[d] \\
N'/\Z \mathbf{1} \arrow[r] & N/\Z \mathbf{a} \arrow[r] & 
N''/\Z \mathbf{1}.
\end{tikzcd}
\end{equation}
Here, $\mathbf{1} = \begin{pmatrix}
    1 \\ 1 \\ \vdots \\ 1 \end{pmatrix}$ is the vector of $1$'s and 
    $\mathbf{a} = \begin{pmatrix}
    a_0 \\ a_1 \\ \vdots \\ a_{n+1} \end{pmatrix}$ is the vector
of weights.  The lattice maps from the top row descend to the bottom
because
$$B \mathbf{1} = dA^{-1} \mathbf{1} = \frac{d}{m}\mathbf{a},$$
while
$$A \mathbf{a} = A (m A^{-1} \mathbf{1}) = m \mathbf{1}.$$
Finally, the rational
maps of weighted projective spaces in the third row 
of \eqref{eq:toricsetup} also induce
rational maps $\varphi: Y \dashrightarrow X$ and $\psi: X \dashrightarrow Z$
because the corresponding ring homomorphisms 
$\C[M''] \rightarrow \C[M] \rightarrow \C[M']$ send
$z_0 + \cdots + z_{n+1} \mapsto f \mapsto y_0^d + \cdots + y_{n+1}^d$.
We've used here that the composition of these homomorphisms is the 
``$d$th power map"
$z_i \mapsto y_i^d$ since 
$B^{\mathsf{T}}A^{\mathsf{T}} = (AB)^{\mathsf{T}} = dI_{n+2}$.

When there is no possibility for confusion, we will write expressions
such as $N''/N, M'/M$, etc., instead of $N''/AN, M'/B^{\mathsf{T}}M$ to
denote the quotient of one of these lattices by the sublattice of finite
index which is the image of one of the maps above.

Each horizontal map in the diagram \eqref{eq:toricsetup} is 
(birational to) a quotient by a finite subgroup of the torus action. 
Specifically, this subgroup is given by the quotient
of the corresponding cocharacter lattices \cite[p. 34]{Fulton}.
Therefore, the rational map $\varphi$ in \eqref{eq:toricsetup} is 
the quotient of $X$ by the finite
abelian group 
\begin{equation}
\label{eq:Gamma_def}
\Gamma \coloneqq (N/\Z \mathbf{a})/B(N'/\Z \mathbf{1}).
\end{equation}
We may identify $\Gamma$ with a subgroup of 
$\mathrm{Aut}_{\mathrm{tor}}(Y) \cong (\Z/d)^{n+2}/(\Z/d)$
via
\begin{align*}
 (N/\Z \mathbf{a})/B(N'/\Z \mathbf{1}) & \cong 
A(N/\Z \mathbf{a})/(AB(N'/\Z \mathbf{1})) \\
 & = A(N/\Z \mathbf{a})/(dN''/\Z \mathbf{1}) \subseteq 
(N''/\Z \mathbf{1})/(dN''/\Z \mathbf{1}) \cong (\Z/d)^{n+2}/(\Z/d).   
\end{align*}

The map of tori $\Spec(\C[M]) \rightarrow \Spec(\C[M''])$ inducing
the map $\psi$ in \eqref{eq:toricsetup} is naturally identified
as a quotient map with respect
to the action of the finite abelian group $N''/N$. Viewing
$N''$ as the superlattice $A^{-1} \Z^{n+2}$ of $N$, we see that this is
exactly the quotient by 
$A^{-1} \Z^{n+2}/\Z^{n+2} \cong \mathrm{Aut}_{\mathrm{tor}}(f)$.
On the level of weighted projective spaces, we see that
$\psi: X \dashrightarrow Z$ is birational to the quotient by the 
full group $\mathrm{Aut}_{\mathrm{tor}}(X)$.  By exactly the same
reasoning, $\psi \circ \varphi: Y \dashrightarrow Z$ is the 
quotient by $\mathrm{Aut}_{\mathrm{tor}}(Y) \cong (\Z/d)^{n+2}/(\Z/d)$.

We summarize the results of this subsection in the following lemma.

\begin{lemma}
\label{lem:quotient}
Suppose that $X$ is an irreducible
Delsarte hypersurface of
dimension $n$.
Then there is a 
natural sequence of dominant rational maps
\begin{equation*}
\begin{tikzcd}
Y \arrow[r,dashed,"\varphi"] & X \arrow[r,dashed,"\psi"] & Z     
\end{tikzcd}
\end{equation*}
where $Y \subseteq \mathbb{P}^{n+1}$ is the Fermat hypersurface of degree 
$d$ and $Z \subseteq \mathbb{P}^{n+1}$ is the Fermat hypersurface (hyperplane)
of degree
$1$.  The map $\varphi$ is birational to the quotient of the
Fermat surface by the group $\Gamma
\subset \mathrm{Aut}_{\mathrm{tor}}(Y)$ defined in \eqref{eq:Gamma_def}. 
Each of the maps $\psi: X \dashrightarrow Z$ and
$\psi \circ \varphi: Y \dashrightarrow Z$ is birational to the quotient
of the respective
Delsarte hypersurface by its full toric automorphism group.
\end{lemma}

\subsection{Cohomology and birational invariants}
\label{subsect:cohom}
For an irreducible variety $V$, we make the following definitions,
following Shioda:

$$H^n(k(V)) \coloneqq \varinjlim_{U \subseteq V \mathrm{ open}} H^n(U),
\hspace{1cm} T^n(V) = \im(H^n(V) \rightarrow H^n(k(V))).$$
Here $H^n$ stands
for de Rham cohomology $H^n(V,\C)$ with complex coefficients.
In modern terminology, $T^n(V)$ is the first subquotient $N^0 H^n/N^1 H^n$ of
the coniveau filtration on the cohomology group $H^n$. 

The groups $H^n(k(V))$ and $T^n(V)$ are birational invariants of irreducible
varieties.  Also, $T^n(V)$ is well-behaved under finite quotients:

\begin{proposition}[{{{{\cite[Proposition 5]{Shioda86}}}}}]
\label{prop:Transcendental_inv}
Let $\Gamma$ be a finite group of automorphisms of $V$.  Then 
$T^n(V/\Gamma) = T^n(V)^{\Gamma}$.
\end{proposition}

For a normal variety $V$ and a positive integer $m$,
the \textit{$m$th plurigenus} of $V$ is 
$H^0(V,mK_V)$, where $K_V$ is the 
canonical divisor of $V$.  The 
\textit{Kodaira dimension} of $V$ is
the minimum $\kappa = \kappa(V)$
such that $h^0(V,mK_V)/m^{\kappa}$ is bounded as
$m \rightarrow \infty$
(we set $\kappa(V) = -\infty$ if all plurigenera are $0$).
The plurigenera are birational invariants for varieties with
canonical singularities, but only the first, also called the
\textit{geometric genus} $p_g(V) \coloneqq h^0(V,K_V)$, is a 
birational invariant of klt varieties.
The \textit{canonical map} of $V$ is the rational map
$V \dashrightarrow \mathbb{P}H^0(V,K_V)$.  This is also
a birational invariant of klt varieties.

When $V$ is a smooth surface, 
$\mathrm{rank}(T^2(V)) = b_2 - \rho = \lambda$, where 
$b_2$ is the second Betti number of $V$, 
$\rho$ is the Picard number $\mathrm{rk}(NS(V))$, and 
their difference $\lambda$ is
called the \textit{Lefschetz number} \cite[Remark 6]{Shioda86}.  The inequality
$\rho \leq h^{1,1}$ and the Hodge decomposition give
$$\lambda = b_2 - \rho \geq b_2 - h^{1,1} = 2h^{2,0}.$$
A smooth surface has \textit{maximal Picard number} when $\rho = h^{1,1}$, which
is equivalent to $\lambda = 2h^{2,0} = 2 h^0(V,K_V) = 2p_g$ by the above.
We say that an irreducible
surface has \textit{maximal Picard number} when a resolution of singularities
(equivalently every resolution of singularities) does so.

\begin{remark}
\label{rem:pg0-maxPicardNumber}
If $X$ is a smooth surface with geometric genus $p_{g} = 0$, then
$X$ automatically has maximal Picard number.  Indeed, we have
$H^{2,0}(X) = H^0(X,K_X) = 0$ so
the Hodge decomposition gives $H^{2}(X, \mathbb{C}) = H^{1, 1}(X)$.
By the Lefschetz (1,1)-theorem, $\Pic(X) \rightarrow H^{1, 1}(X, Z)$ 
is surjective, so $H^{1,1}(X)$ is spanned, as a $\C$-vector space,
by the classes of divisors on $X$. 
This shows $H^{1,1}(X) \cong NS(X) \otimes \mathbb{C}$, and therefore
$\rho(X) = h^{1,1}(X)$.
\end{remark}

In this paper, we consider normal surfaces with quotient singularities.
Note that a normal surface singularity is klt if and only 
it is a quotient singularity \cite[Theorem 7.4.11]{Ishii}.
The criterion $\lambda = 2p_g$ 
for having maximal Picard number 
will be convenient to check on klt models:

\begin{lemma}
\label{lem:max_criterion}
Suppose $V$ is a surface with quotient singularities.  
Then $V$ has maximal Picard
number if and only if $\lambda(V) = 2p_g(V)$.
\end{lemma}

\begin{proof}
The space $H^0(V,K_V)$ is a birational invariant of varieties with klt 
singularities, while $\lambda(V)$ is a birational invariant of 
irreducible varieties, so both these invariants equal their values on a 
nonsingular model $\widetilde{V}$ of $V$.  By definition, $V$ has maximal
Picard number if and only if $\widetilde{V}$ does, which is equivalent to
$\lambda(V) = \lambda(\widetilde{V}) = 2p_g(\widetilde{V}) = 2p_g(V)$.
\end{proof}

Finally, we note that
the space of sections of the canonical bundle is well-behaved under
finite quotients:

\begin{lemma}
\label{lem:quotient_forms}
Let $V$ be a smooth variety, $\Gamma$ a finite group of automorphisms of
$V$, and $W \rightarrow V/\Gamma$ a resolution of singularities.  Then
$H^0(W,K_W) \cong H^0(V,K_V)^{\Gamma}$.
\end{lemma}

\begin{proof}
Let $U \subseteq V/\Gamma$ be the smooth locus of $V/\Gamma$ and $\pi^{-1}(U)$
its preimage in $V$ under the quotient morphism $\pi: V \rightarrow V/\Gamma$.
Both $\pi^{-1}(U)$ and $U$ have complement of codimension at least $2$ in $V$
and $V/\Gamma$, respectively, so it's enough to show that 
$H^0(U,K_U) \cong H^0(\pi^{-1}(U),K_{\pi^{-1}(U)})^{\Gamma}$.

The restricted quotient 
morphism $\pi|_{\pi^{-1}(U)}: \pi^{-1}(U) \rightarrow U$ is a finite map of
smooth varieties, which gives rise to a pullback morphism of sheaves of 
differential forms $\pi^*: \omega_U \rightarrow \pi_*(\omega_{\pi^{-1}(U)})$.
Since we're working in characteristic zero, the morphism $\pi^*$ induces
an isomorphism $\omega_U \cong (\pi_*(\omega_{\pi^{-1}(U)}))^\Gamma$
(see, e.g., \cite[Theorem 1]{Brion}). Taking 
spaces of global sections gives the desired result.

\end{proof}

\section{Computing the Lefschetz number}
\label{sect:algorithm}

In this section, we lay out the algorithm for computing 
the Lefschetz number of a Delsarte surface in weighted
projective space.  The exposition parallels
the original work of Shioda \cite{Shioda82, Shioda86}, which
in turn relies on the theory of Poincar\'{e} residues
(see more on the history below).
However, we present the results in the toric language of 
\Cref{sect:prelim}, which ultimately leads to some new insights.

\subsection{The cohomology of the Fermat hypersurface}
\label{subsect:Fermat}

Let $Y$ be the Fermat hypersurface of dimension $n$ and degree $d$.
The large symmetry group of $Y$ adds enough structure to the middle
cohomology $H^n(Y,\C)$ that the space of Hodge classes on $Y$ may be 
precisely computed.  The first step is to describe the action of 
the toric automorphism group 
$\mathrm{Aut}_{\mathrm{tor}}(Y) = 
(\Z/d)^{n+2}/(\Z/d)$ on the cohomology group $H^n(Y,\C)$.  

To do this, we use the description of the cohomology of 
hypersurfaces in projective space using Poincar\'{e} residues, which is
due to Griffiths \cite{Griffiths}.
For any smooth hypersurface $Y \subseteq \mathbb{P}^{n+1}$ of dimension 
$n$ defined by an equation $f$, and every $0 \leq p \leq n$, there is an
isomorphism
$$R^{t(p)} \xrightarrow{\cong} H_{\mathrm{prim}}^{n-p,p}(Y),$$
where $R$ is the quotient of the polynomial ring in $n+2$ variables by the 
Jacobian ideal $J$ generated by partial derivatives of $f$,
and $R^{t(p)}$ is the graded piece in degree
$t(p) \coloneqq (p+1)d - (n+2)$.  The map 
takes (the equivalence class of) a polynomial $\alpha$ and sends it to the 
Poincar\'{e} residue
of the rational $(n+1)$-form $\Omega_{\alpha}$
on projective space $\mathbb{P}^{n+1}$ given by
$$\Omega_{\alpha} \coloneqq \frac{\alpha \Omega}{f^{p+1}},
\text{ where }\Omega = \sum_{i = 0}^{n+1} (-1)^i 
y_i dy_0 \wedge \cdots \wedge \widehat{dy_i} \wedge \cdots \wedge dy_{n+1}.$$
We may readily compute the action of $\mathrm{Aut}(Y)$ on 
polynomials of a given degree (as well as the action on $\Omega$),
which yields a description of the $\mathrm{Aut}(Y)$-action 
on $H_{\mathrm{prim}}^{n-p,p}(Y)$.

In the special case that $Y$ is the Fermat hypersurface
$Y = \{y_0^d + \cdots + y_{n+1}^d = 0\}$, we have
$J = (y_0^{d-1},\ldots,y_{n+1}^{d-1})$, so a natural basis
of $R^{t(p)}$ consists of
the monomials $\alpha = y_0^{\alpha_0-1} \cdots y_{n+1}^{\alpha_{n+1}-1}$,
where $\alpha_0,\ldots,\alpha_{n+1} \in \{1,\ldots,d-1\}$ and 
$\alpha_0 + \cdots + \alpha_{n+1} = (p+1)d$.  We'll use the notation $|\alpha|$
for the sum of coordinates of the tuple $(\alpha_0,\ldots,\alpha_{n+1})$.

With the same setup and notation as \Cref{subsect:quotient_setup}, 
characters of the toric automorphism
group $N''/N' \cong (\Z/d)^{n+2}$ of the Fermat equation $f_{F}$
are naturally viewed as 
elements of the quotient of dual lattices $M'/M''$, 
which we may think of as the group of integer
lattice points in the cube $[0,d]^{n+2}$ in $M'$, with addition taken modulo $d$.
The character group of the toric automorphism group 
$\mathrm{Aut}_{\mathrm{tor}}(Y)$ (which is a quotient of 
$\mathrm{Aut}_{\mathrm{tor}}(f_{F})$) is the subgroup of 
$M'/M''$ consisting of those lattice points whose coordinates
sum to a multiple of $d$.  From now on, we identify characters and lattice
points in this way.

When $\alpha = y_0^{\alpha_0-1} \cdots y_{n+1}^{\alpha_{n+1}-1}$,
the cohomology class $\mathrm{Res}(\Omega_{\alpha}) \in H^n(Y,\C)$ is an 
eigenvector for the action of $\mathrm{Aut}_{\mathrm{tor}}(Y)$.  The 
character for this eigenvector corresponds to the
lattice point $(\alpha_0,\ldots,\alpha_{n+1})$ under the identification in the
last paragraph,
because $\mathrm{Aut}_{\mathrm{tor}}(Y)$ acts on $\Omega$ as the 
determinantal character, which we can think of as $y_0 \cdots y_{n+1}$.
The classes $\mathrm{Res}(\Omega_{\alpha})$ form a basis for
$H^n(X,\C)$ as $\alpha$ ranges over tuples $(\alpha_0,\ldots,\alpha_{n+1})$ 
with $|\alpha| \equiv 0 \pmod d$.

Since Hodge numbers are deformation invariants in smooth families, this 
also gives a nice lattice-theoretic interpretation of the Hodge numbers of any
smooth hypersurface of degree $d$ in $\mathbb{P}^{n+1}$:
to find $h^{p,q}_{\mathrm{prim}}$,
we simply add up the number of lattice points in the interior of the
cube $[0,d]^{n+2}$
that lie on the slice where the sum of coordinates is $(p+1)d$.

\begin{example}
Let $X_{4}$ be a smooth hypersurface of degree $d = 4$ in
$\mathbb{P}^{3}$ $(n = 2)$. Then, there are two hyperplanes of interest:

\begin{equation*}
    \begin{cases}
     \alpha_0 + \alpha_1 + \alpha_2 = 4 \\
     \alpha_0 + \alpha_1 + \alpha_2 = 8 \\
    \end{cases}   
\end{equation*}

The figure below shows these planes in the cube $[0, 4]^{3}$,
where $\alpha_0 + \alpha_1 + \alpha_2 = 4$ is represented by 
the blue plane passing through points $(4, 0, 0)$, $(0, 4, 0)$,
and $(0, 0, 4)$, and $\alpha_0 + \alpha_1 + \alpha_2 = 8$ is 
represented by the red plane passing through $(4, 4, 0)$, 
$(0, 4, 4)$, and $(4, 0, 4)$.

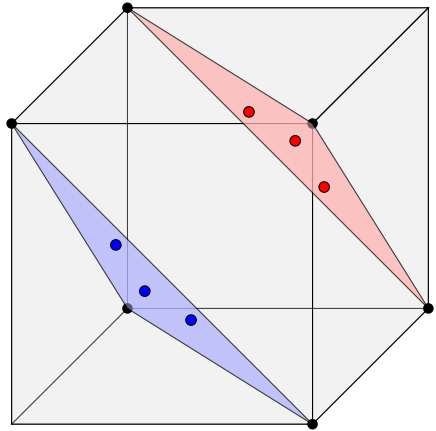
\begin{figure}[h]
\begin{tikzpicture} \label{fig:lattice-d2-n1}
    \draw[fill=gray!20, opacity=0.5] (0,0,-4) -- (4,0,-4) -- (4,4,-4) -- (0,4,-4) -- cycle; 
    \draw[fill=gray!20, opacity=0.5] (0,0,0) -- (0,0,-4) -- (0,4,-4) -- (0,4,0) -- cycle;   
    \draw[fill=gray!20, opacity=0.5] (0,0,0) -- (0,0,-4) -- (4,0,-4) -- (4,0,0) -- cycle;   

    \draw (4,0,0) -- (4,0,-4) -- (4,4,-4) -- (4,4,0) -- cycle; 
    \draw (0,4,0) -- (4,4,0) -- (4,4,-4) -- (0,4,-4) -- cycle; 
    \draw (0,0,0) -- (4,0,0) -- (4,4,0) -- (0,4,0) -- cycle;   

    \draw[fill=blue!30, opacity=0.5] (4,0,0) -- (0,0,-4) -- (0,4,0) -- cycle;

    \draw[fill=red!30, opacity=0.5] (4,4,0) -- (4,0,-4) -- (0,4,-4) -- cycle;

    \fill[black] (4,0,0) circle (2pt);
    \fill[black] (0,0,-4) circle (2pt);
    \fill[black] (0,4,0) circle (2pt);
    \fill[black] (4,4,0) circle (2pt);
    \fill[black] (4,0,-4) circle (2pt);
    \fill[black] (0,4,-4) circle (2pt);

    \coordinate (C1) at (0,0,0);
    \coordinate (C2) at (4,0,0);
    \coordinate (C3) at (0,4,0);
    \coordinate (C4) at (0,0,-4);
    \draw [fill=blue!30, opacity=0.5] (C2) -- (C3) -- (C4) -- cycle;
    \draw[fill = blue, draw = black] (1,1,-2) circle (2pt);
    
    \draw[fill = blue, draw = black] (2,1,-1) circle (2pt);
    \draw[fill = blue, draw = black] (1,2,-1) circle (2pt);

    \coordinate (C1) at (0,0,0);
    \coordinate (C2) at (4,0,-4);
    \coordinate (C3) at (4,4,0);
    \coordinate (C4) at (0,4,-4);
    \draw [fill=red!30, opacity=0.5] (C2) -- (C3) -- (C4) -- cycle;
    \draw[fill = red, draw = black] (3,2,-3) circle (2pt);
    \draw[fill = red, draw = black] (2,3,-3) circle (2pt);

    \draw[fill = red, draw = black] (3,3,-2) circle (2pt);
\end{tikzpicture}
\end{figure}

The number of integer lattice points on the blue plane in the
interior of this cube gives us the value of $h^{1,0}_{\rm prim}(X_{4})$,
and the number of interior integer lattice points on the red plane gives
us $h^{0,1}_{\rm prim}(X_{4})$. In this case, we have interior points
$(1,1,2), (2,1,1), (1,2,1)$ on the blue plane, and interior points
$(3,2,3), (2,3,3), (3,3,2)$ on the red plane, so that
\begin{center}
$h^{1,0}(X_4) = h^{1,0}_{\rm prim}(X_{4}) = 3 
= h^{0,1}_{\rm prim}(X_{4}) = h^{0,1}(X_4)$.
\end{center}
\end{example}

To determine the Hodge classes on the Fermat hypersurface,
we need to 
additionally understand the action of $\Gal(\C/\Q)$ on the direct
sum decomposition of $H^n(Y,\C)$ into eigenspaces.  
The eigenspace decomposition is already defined
over $\Q(\zeta_d)$, where $\zeta_d$ is a primitive $d$th root of unity,
so in fact we can reduce to understanding the action
of the Galois group $\Gal(\Q(\zeta_d)/\Q) = (\Z/d)^*$.
An element $t \in (\Z/d)^*$ maps $\xi$ to $\xi^t$ for each 
$d$th root of unity $\xi$. 
There is an induced action on the group of characters by 
mapping $\alpha \mapsto t \alpha$, where 
$t \cdot (\alpha_0,\ldots,\alpha_{n+1}) = 
(t\alpha_0 \bmod d,\ldots, t\alpha_{n+1} \bmod d)$.
When the dimension $n = 2p$ is even, the space of Hodge classes 
(tensor $\C$) is a sum of those eigenspaces for which $|t \alpha| = (p+1)d$
for every $t \in (\Z/d)^*$.  All these results are summarized in the 
following theorem, which is stated in \cite[Theorem I]{Shioda79},
though parts of it date back to earlier work 
(see \cite[Section 6]{Katz} or \cite[Section 3]{Ogus}).
Recall in the theorem below that the notation $|\alpha|$ refers
to the sum of the coordinates of a tuple $(\alpha_0,\ldots,\alpha_{n+1})$,
where all coordinates are between $1$ and $d-1$.

\begin{theorem}
\label{thm:cohomdecomp}
Let $Y$ be the Fermat hypersurface of degree $d$ in $\mathbb{P}^{n+1}$.
Then the action of $\mathrm{Aut}_{\mathrm{tor}}(Y) \cong (\Z/d)^{n+2}/(\Z/d)$
on $H^n_{\mathrm{prim}}(Y,\C)$ decomposes as
$$H^n_{\mathrm{prim}}(Y,\C) = \bigoplus_{\alpha} V(\alpha)$$
where each $V(\alpha)$ is a one-dimensional eigenspace acted on by
the corresponding character, as $\alpha$ ranges over all tuples
$(\alpha_0,\ldots,\alpha_{n+1})$ where
$\alpha_0,\ldots,\alpha_{n+1} \in \{1,\ldots,d-1\}$ and
$d$ divides $|\alpha|$.
For each $p$,
$H_{\mathrm{prim}}^{p,n-p}(Y)$ is the direct sum of those
summands $V(\alpha)$ with
$|\alpha| = (p+1)d$.
Further, when $n = 2 \ell$, we have that the space of Hodge classes
$$(H_{\mathrm{prim}}^n(Y,\Q) \cap H^{\ell,\ell}(Y)) \otimes_{\Q} \C$$
is the sum of the summands $V(\alpha)$ for which 
$|t \alpha| = (\ell+1)d$
for all $t \in (\Z/d)^*$.
\end{theorem}

\begin{example}
Let $X_5$ be a smooth hypersurface of degree $d = 5$ in
$\mathbb{P}^{3}$. Then, there are three interior
``slices" of the cube
$[0,5]^4$ where the coordinates add up to a multiple of $5$:

\begin{equation*}
    \begin{cases}
      \alpha_0 + \alpha_1 + \alpha_2 + \alpha_3 = 5 \\
      \alpha_0 + \alpha_1 + \alpha_2 + \alpha_3 = 10  \\
      \alpha_0 + \alpha_1 + \alpha_2 + \alpha_3 = 15 
    \end{cases}   
\end{equation*}

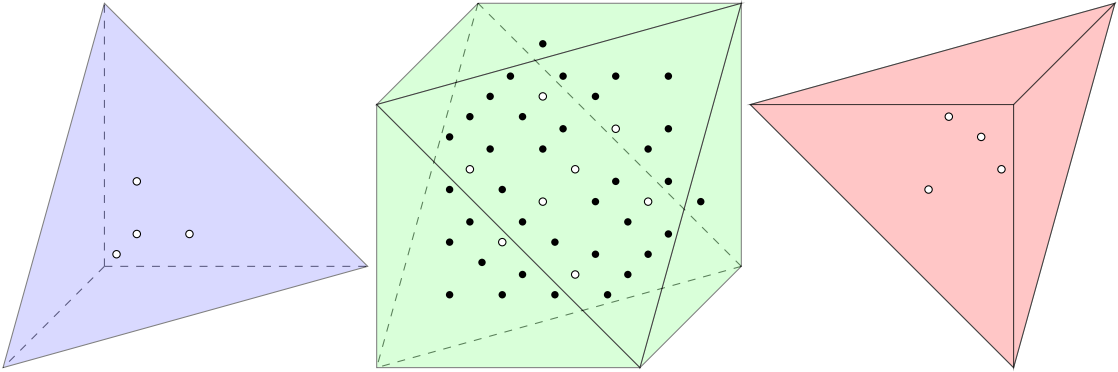
\begin{figure}[h]
\begin{tikzpicture}[scale = 0.7]
\coordinate (A1) at (0,0,0);
\coordinate (A2) at (5,0,0);
\coordinate (A3) at (0,5,0);
\coordinate (A4) at (0,0,5);
\draw [dashed] (A1) -- (A2);
\draw [dashed] (A1) -- (A3);
\draw [dashed] (A1) -- (A4);
\draw [fill=blue!30, opacity=0.5] (A2) -- (A3) -- (A4) -- cycle;
\draw[fill = white, draw = black] (1,1,1) circle (2pt);
\draw[fill = white, draw = black] (2,1,1) circle (2pt);
\draw[fill = white, draw = black] (1,2,1) circle (2pt);
\draw[fill = white, draw = black] (1,1,2) circle (2pt);

\end{tikzpicture}
\begin{tikzpicture}[scale = 0.7]
\coordinate (B1) at (5,0,0);
\coordinate (B2) at (0,5,0);
\coordinate (B3) at (0,0,5);
\coordinate (B4) at (5,5,0);
\coordinate (B5) at (5,0,5);
\coordinate (B6) at (0,5,5);
\draw [dashed] (B1) -- (B2);
\draw [dashed] (B1) -- (B3);
\draw [dashed] (B2) -- (B3);
\draw [fill=green!30, opacity=0.5] (B3) -- (B5) -- (B6) -- cycle;
\draw [fill=green!30, opacity=0.5] (B2) -- (B4) -- (B6) -- cycle;
\draw [fill=green!30, opacity=0.5] (B1) -- (B4) -- (B5) -- cycle;
\draw [fill=green!30, opacity=0.5] (B4) -- (B5) -- (B6) -- cycle;
\fill[black] (1,1,4) circle (2pt);
\fill[black] (1,4,1) circle (2pt);
\fill[black] (4,1,1) circle (2pt);

\fill[black] (1,2,3) circle (2pt);
\fill[black] (1,3,2) circle (2pt);
\fill[black] (2,1,3) circle (2pt);
\fill[black] (3,1,2) circle (2pt);
\fill[black] (2,3,1) circle (2pt);
\fill[black] (3,2,1) circle (2pt);

\fill[black] (1,2,4) circle (2pt);
\fill[black] (1,4,2) circle (2pt);
\fill[black] (2,1,4) circle (2pt);
\fill[black] (4,1,2) circle (2pt);
\fill[black] (2,4,1) circle (2pt);
\fill[black] (4,2,1) circle (2pt);

\draw[fill = white, draw = black] (1,3,3) circle (2pt);
\draw[fill = white, draw = black] (3,1,3) circle (2pt);
\draw[fill = white, draw = black] (3,3,1) circle (2pt);

\fill[black] (3,2,2) circle (2pt);
\fill[black] (2,3,2) circle (2pt);
\fill[black] (2,2,3) circle (2pt);

\fill[black] (1,3,4) circle (2pt);
\fill[black] (1,4,3) circle (2pt);
\fill[black] (3,1,4) circle (2pt);
\fill[black] (4,1,3) circle (2pt);
\fill[black] (3,4,1) circle (2pt);
\fill[black] (4,3,1) circle (2pt);

\draw[fill = white, draw = black] (2,2,4) circle (2pt);
\draw[fill = white, draw = black] (2,4,2) circle (2pt);
\draw[fill = white, draw = black] (4,2,2) circle (2pt);

\fill[black] (2,3,4) circle (2pt);
\fill[black] (2,4,3) circle (2pt);
\fill[black] (3,2,4) circle (2pt);
\fill[black] (4,2,3) circle (2pt);
\fill[black] (3,4,2) circle (2pt);
\fill[black] (4,3,2) circle (2pt);

\fill[black] (1,4,4) circle (2pt);
\fill[black] (4,1,4) circle (2pt);
\fill[black] (4,4,1) circle (2pt);

\fill[black] (2,2,5) circle (2pt);
\fill[black] (2,5,2) circle (2pt);
\fill[black] (5,2,2) circle (2pt);

\draw[fill = white, draw = black] (2,2,2) circle (2pt);

\draw[fill = white, draw = black] (3,3,3) circle (2pt);

\end{tikzpicture}
\begin{tikzpicture}[scale = 0.7]
\coordinate (C1) at (5,5,0);
\coordinate (C2) at (5,0,5);
\coordinate (C3) at (0,5,5);
\coordinate (C4) at (5,5,5);
\draw [fill=red!30, opacity=0.5] (C1) -- (C2) -- (C3) -- cycle;
\draw [fill=red!30, opacity=0.5] (C1) -- (C2) -- (C4) -- cycle;
\draw [fill=red!30, opacity=0.5] (C1) -- (C3) -- (C4) -- cycle;
\draw [fill=red!30, opacity=0.5] (C2) -- (C3) -- (C4) -- cycle;
\draw[fill = white, draw = black] (4,3,3) circle (2pt);
\draw[fill = white, draw = black] (3,4,3) circle (2pt);
\draw[fill = white, draw = black] (3,3,4) circle (2pt);
\draw[fill = white, draw = black] (4,4,4) circle (2pt);
\end{tikzpicture}
\end{figure}

The diagram shows these three slices
(projected to lie in 
three dimensions) along with their interior lattice points.
Counting the points on each slice gives the primitive 
Hodge numbers of any smooth quintic surface:
$h^{2,0}(X_5) = 4 = h^{0,2}(X_5)$ and 
$h^{1,1}_{\mathrm{prim}}(X_5) = h^{1,1}(X_5) - 1 = 44$.

The Galois action on these lattice points also determines
the Picard rank of the Fermat quintic $X_5$.
Points in white above are those that move
between the middle and outer 
slices under the Galois action, while points
in black always remain on the middle slice.  There
are $36$ black points and $16$ white points so
$\rho(X_5) = 36 + 1 = 37$ and $\lambda(X_5) = 16$.
\end{example}

\subsection{The algorithm}
\label{subsect:algorithm}

For simplicity, we now specialize to the situation of surfaces,
with $n = 2$.  Following the setup of \Cref{sect:prelim}, 
suppose we are given a Delsarte equation $f$ defining a 
(well-formed) Delsarte surface $X \subseteq \mathbb{P}(a_0,a_1,a_2,a_3)$
with quotient singularities.
We continue to use the notation from \Cref{subsect:quotient_setup},
so in particular, $X$ is birational to a 
quotient $Y/\Gamma$ of a Fermat surface of degree $d$.
The idea of the formula for the Lefschetz number of $X$ is that
the results from \Cref{subsect:Fermat}
carry over from $Y$ to $X$ upon taking $\Gamma$-invariants.

We will use the following notation for subsets of lattice
points related to the ``cube" $M'/M'' \cong (\Z/d)^4$, where
we view $M$ as a sublattice of $M'$:
\begin{align*}
    \mathcal{A}(X) & \coloneqq \{\alpha = 
    (\alpha_0,\alpha_1,\alpha_2,\alpha_3) \in M/M'' : 
    |\alpha| \equiv 0 \bmod d \}; \\
    \mathcal{A}^{\circ}(X) & \coloneqq \{\alpha = 
    (\alpha_0,\alpha_1,\alpha_2,\alpha_3) \in \mathcal{A} : 
    1 \leq \alpha_i \leq d-1, 0 \leq i \leq 3 \}, \\
    S_p(X) & \coloneqq \{\alpha \in \mathcal{A}^{\circ} : |\alpha| = (p+1)d\}.
\end{align*}
We have the following formula for
the Lefschetz number of $X$.

\begin{theorem}
\label{thm:Lefschetz_formula}
Suppose that $X$ is a Delsarte surface with quotient singularities
in weighted projective
$3$-space.  Then the Lefschetz number of $X$ is given by the formula
$$\lambda(X) = |\mathcal{A}^{\circ}| - |\{\alpha: t \alpha \in S_1 
\text{ for all } t \in (\Z/d)^*\}|.$$
In particular, the following conditions are equivalent:
\begin{enumerate}
    \item $X$ has maximal Picard number;
    \item For all $\alpha \in S_1$ and for all $t \in (\Z/d)^*$,
    $t \alpha \in S_1$;
    \item For all $\alpha \in S_0$ and for all $t \in (\Z/d)^*$,
    $t \alpha \in S_0$ or $t \alpha \in S_2$.
\end{enumerate}
\end{theorem}

In summary, this is a completely lattice-theoretic method of computing
the Lefschetz number of $X$, given the matrix associated to its
equation.  We implement this formula in the \texttt{MAGMA} algorithm
\cite{magma}.

\begin{proof}
We saw in \Cref{thm:cohomdecomp} that
$H_{\mathrm{prim}}^{1,1}(Y,\Q) \otimes_{\Q} \C$ is a sum
of eigenspaces $V(\alpha)$ for which 
$\alpha = (\alpha_0,\alpha_1,\alpha_2,\alpha_3)$ is a
lattice point 
in the interior of the cube $[0,d]^4$ and
$|t \alpha| = 2d$ for all 
$t \in (\Z/d)^*$.
The remaining summands $V(\alpha)$
generate the transcendental cohomology
$T^2(Y)$.

Since $X$ has quotient singularities and is 
birational to $Y/\Gamma$, we have
$T^2(X) \cong T^2(Y)^{\Gamma}$
and $H^0(X,K_X) \cong H^0(Y,K_Y)^{\Gamma}$
by \Cref{prop:Transcendental_inv} and \Cref{lem:quotient_forms}, respectively.
Both $T^2(Y)$ and $H^0(Y,K_Y)$ are a sum of some of the $V(\alpha)$,
so we need to figure out which of these summands are
invariant under the subgroup $\Gamma \subseteq \mathrm{Aut}_{\mathrm{tor}}(Y)$.
Identifying each $\alpha$ with a lattice point in $M'$, the invariant
summands are exactly those for which $\alpha$ is contained in 
the sublattice 
$M$.  Indeed, $M \rightarrow M'$ is dual to the 
lattice map $N' \rightarrow N$ that descends to 
$N'/\Z \mathbf{1} \rightarrow N/\Z\mathbf{a}$, where the quotient 
of these last two lattices is $\Gamma$, by 
\eqref{eq:Gamma_def}.  The dual of $N'/\Z \mathbf{1}$
is the group of points of $M'$ with coordinates summing to a multiple
of $d$, while the dual of $N/\Z\mathbf{a}$ can be viewed as the subgroup
of lattice points also belonging to $M$ with the same property.

Thus we've shown that the set $\mathcal{A}^{\circ}$
is a basis for the $\Gamma$-invariant cohomology
$H_{\mathrm{prim}}^2(Y,\C)^{\Gamma}$ of the Fermat surface of degree $d$, 
and $S_p$ is the part of this invariant cohomology 
in the summand $H_{\mathrm{prim}}^{p,q}(Y)$,
i.e., the set of points on the corresponding ``slice".
Now, the transcendental cohomology $T^2(Y)$ has a basis given by
the \textit{complement} 
of the set of all $\alpha$ with nonzero coordinates
in $M'/M''$ and $|t \alpha| = 2d$ for all $t \in (\Z/d)^*$, inside
the set of all $\alpha$ from \Cref{thm:cohomdecomp}.
The Galois action is by multiplication by integers in $M'/M''$, which
must commute with the $\Gamma$-action.  This means that 
$T^2(X) = T^2(Y)^{\Gamma} = T^2(Y) \cap H_{\mathrm{prim}}^2(Y,\C)^{\Gamma}$
has a basis given by the complement of those
points in $\mathcal{A}^{\circ}$ that stay in the middle
slice under the Galois action.  In other words, a basis for $T^2(X)$
is given by the set
$$\mathcal{A}^{\circ} \setminus \{\alpha \in \mathcal{A}^{\circ}: t \alpha \in S_1 
\text{ for all } t \in (\Z/d)^*\}.$$
Taking the order of this set gives the 
formula for the Lefschetz number $\lambda(X)$.

Finally, by \Cref{lem:max_criterion},
$X$ has maximal Picard number if and only if $\lambda(X) = 2p_g(X)$.  We know 
by \Cref{lem:quotient_forms} that 
$p_g(X)$ counts the $\Gamma$-invariant $2$-forms on $Y$, so $p_g(X) = |S_0|$,
hence $2p_g(X) = |S_0| + |S_2|$.  Therefore, the condition
$\lambda = 2p_g$ is equivalent to the statement that
$\{\alpha: t \alpha \in S_1 \text{ for all } 
t \in (\Z/d)^*\} = S_1$, that is, that
the Galois action stabilizes $S_1$.  It's equivalent to
say that
the union of the remaining two slices, $S_0\cup S_2$, is stable
under the Galois action.  Finally, $S_0$ and $S_2$ are exchanged by
conjugation, i.e., the element $-1 \in (\Z/d)^*$, so it's equivalent to
require that the Galois orbit of each point in $S_0$ is contained in
$S_0\cup S_2$.
\end{proof}

\section{Rationality of Geometric Genus Zero Surfaces}
\label{sect:rat}

In this short section, we use results from toric geometry and the
classification of three-dimensional lattice polytopes to prove \Cref{thm:intro_rat}.
Some of the results concerning Newton polytopes will also be used in 
later sections.

A hypersurface in a torus $(\C^*)^{n+1} = 
\Spec(\Z^{n+1}) = \Spec(\C[t_1,t_1^{-1},\ldots,t_{n+1},t_{n+1}^{-1}
])$ is defined by a Laurent polynomial $p$ in the variables $t_1,\ldots,t_{n+1}$. 
The \textit{Newton polytope} of this hypersurface is the convex hull $\Delta$ of 
the points in $\mathbb{Z}^{n+1}$ corresponding to monomials 
with nonzero coefficients in the polynomial $p$.  A polynomial $p$ is
\textit{Newton non-degenerate} if for every face $\delta \subset \Delta$,
the restriction of the polynomial $p$ to $\delta$ defines a smooth hypersurface
of the torus.

Our applications here will deal with polynomials that are general members
of linear systems generated by monomials.  These never have base points
on the torus, so by Bertini's theorem they are always Newton non-degenerate.
For ease of exposition, we'll also assume that all Newton polytopes
are of full dimension, i.e., that $\Delta$ is not contained in any
affine subspace of dimension less than $n+1$.

For a hypersurface $V \coloneqq \{p = 0\} \in (\C^*)^{n+1}$, results of
Khovanskii \cite{Khovanskii83} show that we can always choose an
appropriate projective 
toric variety containing $(\C^*)^{n+1}$ in which the compactification $\bar{V}$
of $V$ is smooth.  Remarkably, we can compute certain birational invariants
of $\bar{V}$ directly from the polytope $\Delta$.

\begin{theorem}[{{{Khovanskii \cite{Khovanskii}}}}]
Let $V \subset (\C^*)^{n+1}$ be a Newton non-degenerate hypersurface with
Newton polytope $\Delta$. The geometric genus $p_g(\bar{V})$ is equal 
to the number of lattice points in the interior of $\Delta$.
\end{theorem}

In the case $n = 1$, this says that the genus of the unique smooth, projective
curve birational to $V$ has genus equal to the number of interior points of
the lattice polygon.  This special case was originally proven by Baker 
\cite{Baker} and is often known as Baker's theorem. We'll make use of it in
\Cref{sect:unbounded_genus}.

We can also measure the Kodaira dimension of $\bar{V}$ using the polytope 
$\Delta$.  This requires a notion called the \textit{Fine interior} of $\Delta$.
The Fine interior $F(\Delta)$ is defined to be the intersection of all half-spaces
$\{v: \langle n, v \rangle  \geq k+1 \}$ such that the half-space
$\{v : \langle n, v \rangle  \geq k\}$
contains
$\Delta$.  Informally, we take all supporting hyperplanes of $\Delta$ and
``slide" them into $\Delta$ until they hit another lattice point.  One can
check that $F(\Delta)$ is a polytope with vertices in $\Q^{n+1}$.

\begin{theorem}[{{{Fine \cite{Fine, Reid85}}}}]
Let $V \subset (\C^*)^{n+1}$ be a Newton non-degenerate hypersurface with
Newton polytope $\Delta$.  Then the Kodaira dimension of $\bar{V}$ is
given by
$$\kappa(\bar{V}) = \begin{cases}
    -\infty, & \textit{if } F(\Delta) = \emptyset, \\
    \dim(F(\Delta)), & \textit{if } 0 \leq \dim(F(\Delta)) \leq n, \\
     \dim(F(\Delta)) - 1, & \textit{if } \dim(F(\Delta)) = n+1.
\end{cases}$$
\end{theorem}

Using these results, we can now prove \Cref{thm:intro_rat}.

\begin{proof}[Proof of \Cref{thm:intro_rat}]
In the above setup, let $V \subset (\C^*)^3$ be a Newton non-degenerate
surface.  A smooth, projective
surface over $\C$ is rational if and only if it has Kodaira dimension
$-\infty$ by the Enriques-Kodaira classification.  Suppose that $V$
satisfies $p_g(\bar{V}) = 0$ but $\kappa(\bar{V}) \geq 0$.  In terms of 
polytopes, this is equivalent to the statement that the three-dimensional
Newton polytope $\Delta$ has no interior 
lattice points (these polytopes are called \textit{hollow}
or \textit{lattice-free} in the references) but has non-empty Fine interior
$F(\Delta)$. We will show that if both conditions are satisfied, 
then the fundamental group of $\bar{V}$ must be nontrivial.

A result of Treutlein \cite{Treutlein} shows that, up to unimodular equivalence,
any three-dimensional
hollow lattice polytope $\Delta$ either (1) projects onto the one-dimensional unit interval,
(2) projects onto the triangle with vertices $(0,0), (2,0), (0,2)$ in $\Z^2$,
or (3) has bounded volume.  If $\Delta$ satisfies (1) or (2), then it has empty
Fine interior.  In case (3), $\Delta$ is contained in some maximal hollow convex
set.  Results of Averkov et al.\ \cite{AWW,AKW} show that all such maximal
sets are lattice polytopes and completely classify them.  It turns out that there 
are exactly $12$ maximal hollow lattice polytopes of dimension $3$.  
In \cite[Section 2.7]{BKS}, the Fine interiors of each are calculated.  Exactly $9$
have nonempty Fine interior, and no proper subpolytopes of any of the $12$ have 
nonempty Fine interior.  In summary, our surface $V$ above must have Newton polytope
which is exactly one of the $9$, up to unimodular equivalence.

Now let $X = X_d \subset \mathbb{P}(a_0,a_1,a_2,a_3)$ be a well-formed surface
with $p_g = 0$.
If any $X$ of degree $d$ in this weighted projective space is quasismooth, the 
general one is; further, the plurigenera of a resolution are the same for any 
quasismooth member.  Hence we may take $X$ to be general so that the hypersurface
it cuts out in the torus $(\C^*)^3$ is Newton non-degenerate.  The Newton polytope
$\Delta$ of $X$ is hollow, but 
cannot be one of the $9$ exceptional examples because all
of these yield surfaces with non-trivial fundamental group \cite[Section 2.7]{BKS},
while the fundamental group of a quasismooth hypersurface of dimension greater 
than $1$ in weighted projective
space is trivial \cite[Theorem 3.2.4]{Dolgachev}.  Thus, $X$ is rational.
\end{proof}

\begin{remark}
A few special cases of \Cref{thm:intro_rat} were known previously.
Given a quasismooth weighted surface $X$,
Chitayat \cite{Chitayat} proved that $p_g = 0$ implies rational 
when all weights divide the degree.
Urz\'ua and Y\'a\~nez \cite{UY} showed the same for Delsarte surfaces
whose equations are ``loops".  These were originally
studied by Koll\'{a}r \cite{Kollar}. 
\end{remark}

For use in the next section, we record the $4$ simplices that appear among
the maximal hollow polytopes of dimension $3$ with nonempty Fine interior
\cite[Table 2.10]{BKS}:

\begin{table}[h]
    \centering
    \begin{tabular}{c|c}
       ID  &  Vertices \\
       \hline
        4 &  $(0,0, 0),(4, 0, 0),(2, 4, 0),(3, 0, 2)$ \\
        \hline
        5 &  $(0, 0, 0),(2, 2, 0),(1, 1, 2),(3,-1, 2)$ \\
        \hline
        9 &  $(0, 0, 0),(3, 0, 0),(1, 3, 0),(2, 0, 3) $ \\
        \hline
        12 & $ (-1,0, 0),(0, 1, -2),(1, 2, 1),(2, -2, -1)$
    \end{tabular}
    \caption{The four lattice simplices in $\mathbb{Z}^3$ with 
    no interior lattice points and nonempty Fine interior.}
    \label{tab:simplices}
\end{table}

\section{Automorphisms and the Lefschetz number}
\label{sect:aut}

We've seen that both the toric automorphism group and the Lefschetz
number
of a Delsarte surface have lattice-theoretic interpretations.  Taken 
together, these descriptions yield a surprising connection between the 
two.  Roughly speaking, a small automorphism group indicates that the
canonical class of a resolution of $X$ is less positive, and $X$ is more
likely to have maximal Picard number, or even be rational.
In this section, we prove \Cref{thm:intro_aut}
and provide a series of examples illustrating the diversity of behavior
exhibited by weighted surfaces with maximal Picard number.

Our first result ties the toric automorphism group 
to the group $\mathcal{A}$ of lattice
points defined in \Cref{subsect:algorithm}.

\begin{lemma}
\label{lem:dual_iso}
Let $X$ be an irreducible
Delsarte surface in weighted projective space.
Then the finite abelian groups $\mathrm{Aut}_{\mathrm{tor}}(X)$ and 
$\mathcal{A}(X)$ are isomorphic.
\end{lemma}

\begin{proof}
We've already seen in \Cref{subsect:quotient_setup}
that the quotient $N''/N$ of 
lattices is isomorphic to $\mathrm{Aut}_{\mathrm{tor}}(f)$, where $f$
is the equation of $X$.  For explicitness, we make an identification 
$N' \cong \Z^4$ so that the sequence
$N' \rightarrow N \rightarrow N''$ is identified with 
$\Z^4 \subseteq B^{-1} \Z^4 \subseteq d^{-1} \Z^4$, 
while the dual sequence $M' \leftarrow M \leftarrow M''$ is identified with
to $\Z^4 \supseteq B^{\mathsf{T}} \Z^4 \supseteq d \Z^4$.

Then we have 
$$N''/N \cong (d^{-1} \Z^4)/B^{-1} \Z^4 = 
B^{-1}A^{-1} \Z^4 / B^{-1} \Z^4 \cong A^{-1} \Z^4/\Z^4 
\cong \mathrm{Aut}_{\mathrm{tor}}(f),$$ 
as expected.  Now, $\mathrm{Aut}_{\mathrm{tor}}(X)$ is the quotient of 
$\mathrm{Aut}_{\mathrm{tor}}(f)$ by the charge vector by 
\Cref{lem:Aut_tor}. In the above coordinates, the charge vector is 
$$v = B^{-1} \frac{1}{m} \begin{pmatrix} a_0 \\ a_1 \\ a_2 \\ a_3 
\end{pmatrix} = \frac{1}{d} \begin{pmatrix} 1 \\ 1 \\ 1 \\ 1
\end{pmatrix},$$
so we could also write $\mathrm{Aut}_{\mathrm{tor}}(X)$ as 
$d^{-1} \Z^4/(B^{-1} \Z^4 + \Z v)$.  This is isomorphic as an
abelian group
to the quotient of the dual of $(B^{-1} \Z^4 + \Z v)$ by the dual
of $d^{-1} \Z^4$. The dual of $B^{-1} \Z^4 + \Z v$ corresponds to 
points $\alpha$ in $B^{\mathsf{T}} \Z^4 \subseteq \Z^4$
that also map $v$ to an integer, i.e., those points of $M$ such that
$\alpha_0 + \alpha_1 + \alpha_2 + \alpha_3$ is a multiple of $d$.
The dual of $d^{-1} \Z^4$, which we identify with $N''$, is $M''$.
Taking the quotient of these gives $\mathcal{A}$, proving the claim.
\end{proof}

When $\mathrm{Aut}_{\mathrm{tor}}(X)$ only has elements of small
order, the same must be true of $\mathcal{A}(X)$, which significantly
constrains the cohomological behavior.

\begin{theorem}
\label{thm:low_exp}
Suppose that $X$ is a Delsarte surface with quotient singularities
in weighted projective $3$-space
and the finite abelian group $\mathrm{Aut}_{\mathrm{tor}}(X)$ has exponent
$e = 1,2,3,4$, or $6$.  Then $X$ has maximal Picard number.  Furthermore,
if $e \leq 3$, then $X$ is rational.
\end{theorem}

\begin{remark}
\begin{enumerate}
\item \Cref{thm:low_exp} generalizes the known fact that the
    Fermat surface $Y_d$ of degree $d = 1,2,3,4$, or $6$ has maximal Picard 
    number, as does any quotient of one of these.
    The theorem applies to many other surfaces only expressible as quotients of
    $Y_d$ for larger $d$, however (see
    the examples later in this section).  
    Conversely, there are also examples where the condition in \Cref{thm:low_exp}
    does not hold but $X$ still has maximal Picard number 
    (see \Cref{sect:unbounded_genus}).
\item  We note that the statement ``$\mathrm{Aut}_{\mathrm{tor}}(X)$ is trivial 
implies $X$ is rational" holds for Delsarte hypersurfaces of 
any dimension. Indeed, the argument for the $e=1$ case below generalizes.
It is a testament to the specialness of Delsarte 
hypersurfaces that a lack of symmetries forces rationality.
The rationality of Delsarte hypersurfaces was previously investigated in 
\cite[Section 3]{Esser} (and in \cite[Section 5]{Kollar} in the special case
of ``loop" hypersurfaces) but these works did not prove the relationship
of rationality to the automorphism group.
\end{enumerate}
\end{remark}

\begin{proof}[Proof of \Cref{thm:low_exp}:]
By \Cref{lem:dual_iso}, 
if $\mathrm{Aut}_{\mathrm{tor}}(X)$ has exponent $e$, the same is true of 
the isomorphic group $\mathcal{A}(X)$. We saw in 
\Cref{thm:Lefschetz_formula} that the Lefschetz number
of $X$ is controlled by the action of the Galois group on the group of
lattice points $\mathcal{A}(X)$.  
Given an element $t \in (\Z/d)^*$, this
action sends any point $\alpha \in \mathcal{A} \subseteq (\Z/d)^4$
to $t \alpha$.  Since the
exponent of $\mathcal{A}$ is $e$ and $e$ divides $d$,
the Galois action by $(\Z/d)^*$ on $\mathcal{A}$
factors through the quotient $(\Z/d)^* \rightarrow (\Z/e)^*$.

Then, if $e = 1,2,3,4$, or $6$, the group of units $(\Z/e)^*$ only has
(at most) the elements $\pm 1$ (i.e., $\phi(e) \leq 2$ where $\phi$ is
Euler totient function).  The action of $1$ is trivial, while the action
of $-1$ is conjugation: it switches points between the sets $S_0$ and $S_2$
but cannot interchange points on either of these slices with 
points on $S_1$.  Thus, \Cref{thm:Lefschetz_formula} shows that $X$
has maximal Picard number.

If $e = 1$, $\mathrm{Aut}_{\mathrm{tor}}(X)$ is trivial.
But it was already established in \Cref{lem:quotient}
that the map $\psi: X \dashrightarrow Z$ is birational to the quotient
of $X$ by $\mathrm{Aut}_{\mathrm{tor}}(X)$, so we conclude that in this 
case $\psi$ is a birational map.  The variety $Z$ is a hyperplane in 
$\mathbb{P}^3$ so it is rational. Hence $X$ is rational.
Note that the requirement that $\mathrm{Aut}_{\mathrm{tor}}(X)$ be trivial
is equivalent to $m = \left|\det(A)\right|$ by \Cref{lem:Aut_tor} and
the fact that $\left|\mathrm{Aut}_{\mathrm{tor}}(f)\right| 
= \left|\det(A)\right|$.
Therefore, this rationality
criterion is equivalent to a result by the first author
\cite[Theorem 3.2]{Esser}.

Finally, we consider the cases $e = 2,3$. 
Since $\mathcal{A}$ must have exponent $e$ also, the options for
lattice points in $\mathcal{A}$ are very limited.
If $e = 2$, then the only lattice point that could be in
$\mathcal{A}^{\circ}$ is
$$\alpha = \begin{pmatrix}
d/2 \\ d/2 \\ d/2 \\ d/2
\end{pmatrix}.$$
This point is in $S_1$, so in particular $S_0$ is empty and
$p_g = 0$.  Similarly, if $e = 3$, then
the only possible points in $\mathcal{A}^{\circ}$ are 
permutations of
$$\alpha = \begin{pmatrix}
d/3 \\ d/3 \\ 2d/3 \\ 2d/3
\end{pmatrix},$$
which are also in $S_1$.  Again, we have $p_g = 0$.

Thus, if $e = 2$ or $3$ and we also assume 
$X$ is quasismooth, \Cref{thm:intro_rat}
proves that $X$ is rational.  We claim that this still holds true
even if we only assume $X$ has quotient singularities.
Indeed, the intersection of $X$ with the torus $(\C^*)^3 \subset \mathbb{P}$
is Newton non-degenerate because a Delsarte surface is always
isomorphic to the general member of the linear system generated by
its monomials.  Because $e = 2$ or $3$ implies $p_g = 0$, we have
that $X$ is rational unless the Newton polytope $\Delta$ of 
$X$ is one of the exceptional polytopes with no interior lattice
points and nonempty Fine interior from \Cref{sect:rat}.  The polytope
$\Delta$ also has $4$ vertices, so it is a simplex, and must therefore be
equivalent to a simplex from \Cref{tab:simplices} under a lattice
isomorphism.

On the other hand, we know that the map $\psi$ from \Cref{lem:quotient}
is the quotient by $\mathrm{Aut}_{\mathrm{tor}}(X)$, so that 
this same group must be the cokernel of the lattice map 
$N/\Z \mathbf{a} \rightarrow N''/\Z \mathbf{1}$ from \eqref{eq:toricsetup_lat}.
The dual of this lattice map takes the unit simplex to $\Delta$
because it maps the equation of a hyperplane to the equation $f$ of $X$.
Hence we can compute what $\mathrm{Aut}_{\mathrm{tor}}(X)$ would have
to be if $\Delta$ were one of the exceptional simplices by taking the quotient
of $\Z^3$ by the sublattice generated by the vectors spanning this simplex. 
A straightforward calculation gives
that these groups for the simplices from \Cref{tab:simplices}
of IDs 4, 5, 9, and 12 are, respectively, $\Z/4 \oplus \Z/8$, $\Z/4 \oplus \Z/4$,
$\Z/3 \oplus \Z/9$, and $\Z/5 \oplus \Z/5$.  None of these has exponent
equal to $2$ or $3$, completing the proof.
\end{proof}

The criterion in \Cref{thm:low_exp} can be applied to a wide range
of surfaces.

\begin{example}
\label{ex:rational_ex}
Consider the surface
$$X_{119} \coloneqq \{x_0^2 x_1 + x_1^3 x_2 + x_2^5 x_3 + x_0 x_3^4 = 0\}
 \subseteq \mathbb{P}(43,33,20,19).$$
This $X$ is quasismooth, well-formed, and has ample
canonical class $K_X = 
\mathcal{O}_X(119 - 43 - 33 - 20 - 19) = \mathcal{O}_X(4)$.
However, a quick matrix computation gives that 
$\mathrm{Aut}_{\mathrm{tor}}(X)$ is trivial, so $X$ is
actually rational (and trivially has maximal Picard number).
\end{example}

\begin{example}
\label{ex:K3_low_exp}
Consider the surface
$$X_{12} \coloneqq \{x_0^3 + x_1^3 + x_2^4 + x_2^2 x_3^6 = 0\}
\subseteq \mathbb{P}(4,4,3,1).$$
One can check that this equation defines a well-formed surface,
which is not quasismooth but still has only quotient singularities.
Since $K_X = 
\mathcal{O}_X(12-4-4-3-1) = \mathcal{O}_X$ is trivial, $X$ is 
in fact a canonical K3 surface.
This $X$ is a quotient of the Fermat surface of degree $12$.
Computing with the associated
matrix gives
$$\mathrm{Aut}_{\mathrm{tor}}(X) \cong \Z/3 \oplus \Z/6,$$
which has exponent $6$, so $X$ has maximal Picard number.
In other words, a minimal resolution of $X$ is a smooth
K3 surface with Picard number $\rho = 20$.
\end{example}

\begin{example}
\label{ex:gen_type_low_exp}
Consider the surface
$$X_{30} \coloneqq \{x_0^5 + x_1^6 + x_2^6 + x_3^6 x_0 = 0\} \subseteq 
\mathbb{P}(6,5,5,4).$$
This surface is well-formed and quasismooth; 
it is a quotient of
the Fermat surface of degree $30$.
A matrix computation gives that 
$$\mathrm{Aut}_{\mathrm{tor}}(X) \cong \Z/6 \oplus \Z/6,$$
which has exponent $6$, so $X$ has maximal Picard number.
The surface $X$ is not canonical at $[0:0:0:1]$, but we can
nevertheless show that a resolution of $X$ is a surface
of general type, i.e., its Kodaira dimension is $2$.
Indeed, the adjunction formula gives 
$K_X = \mathcal{O}_X(30-6-5-5-4) = \mathcal{O}_X(10)$, which
has sections $x_0 x_3, x_1^2, x_1 x_2, x_2^2$, giving $p_g = 4$.
The image of the canonical map $X \dashrightarrow \mathbb{P}^3$
with coordinates $z_0,z_1,z_2,z_3$ corresponding to these sections
is $\mathrm{Proj}(\C[z_0,z_1,z_2,z_3]/(z_1z_3 - z_2^2))$, which we
can identify with the projective 
cone over a conic curve in $\mathbb{P}^2$.
Since the canonical map has image of dimension $2$ and
this image is a birational invariant of klt varieties, a 
resolution of $X$ is of general type.

In fact, we claim that a minimal model of $X$ is a Horikawa surface
with $p_g = 4, K^2 = 4$.  Indeed, $X$ has a unique non-canonical
singularity at the point $p = [0:0:0:1]$, which is of type $\frac{1}{4}(1,1)$.
This singularity is resolved by a single blowup of $p$ with 
exceptional locus a curve $E$ of self-intersection $-4$.
A minimal resolution $\pi: \widetilde{X} \rightarrow X$ of $X$ satisfies
$$K_{\widetilde{X}} = \pi^* K_X - \tfrac{1}{2}E,$$
so that 
$$K_{\widetilde{X}}^2 = K_X^2 + \tfrac{1}{4}E^2 = 
\frac{30 \cdot 10^2}{6 \cdot 5 \cdot 5 \cdot 4} + \tfrac{1}{4}(-4) = 4.$$
It's not hard to see that $\widetilde{X}$ is a minimal surface.
This $\widetilde{X}$ is on the Noether line $p_g = \tfrac{1}{2}K^2 + 2$,
so it is a Horikawa surface.  Horikawa \cite{Horikawa} determined
the possible canonical maps for such surfaces, and the canonical 
images are always surfaces of degree $p_g - 2$ in $\mathbb{P}^{p_g-1}$.
In this case, we saw above that we get a degree $2$ surface in $\mathbb{P}^3$
as the canonical image.
Horikawa surfaces with
maximal Picard number 
were known to exist by the work of Persson \cite{Persson}.
\end{example}

\section{Examples with unbounded geometric genus}
\label{sect:unbounded_genus}

Results of Heijne \cite{Heijne16} show that
Delsarte surfaces with ADE singularities in
$\mathbb{P}^3$ with maximal Picard number all have degree at most $6$, and in
particular have $p_g \leq 10$. 
Among weighted Delsarte surfaces, it's not hard to use \Cref{thm:low_exp}
to find examples of arbitrarily large degree, or arbitrarily large
$d/\max\{a_0,\ldots,a_{n+1}\}$, which have $p_g = 0$.  These
of course have maximal
Picard number.  

However, we can find more interesting behavior as well.
In fact, there exist weighted Delsarte surfaces
with arbitrarily large geometric genus achieving maximal Picard 
number.
More strongly, we find an example for every possible value of $p_g$;
we can even arrange that it is quasismooth.

\begin{theorem}
\label{thm:elliptic_ex}
For every positive integer $s$, there exists a
quasismooth Delsarte 
surface $X = X_s$
in weighted projective $3$-space with geometric genus $p_g = s$ 
and maximal Picard number such that a resolution of
$X_s$ is a simply connected
elliptic surface.
\end{theorem}

Note that for $s \geq 2$, $K_{X_s}$ is ample, but the 
desingularization of $X_s$ is an elliptic surface,
with Kodaira dimension 
$\kappa = 1$.
This illustrates how non-canonical Delsarte hypersurfaces can
achieve intermediate Kodaira dimension, unlike in the canonical case.
Computing the Kodaira dimension of a non-canonical weighted hypersurface
is more difficult; indeed, this is the hardest part of the proof.
Persson \cite[Section 5.5]{Persson}
previously constructed elliptic surfaces of 
arbitrary $p_g$ with maximal Picard number as double covers of 
$\mathbb{P}^1 \times \mathbb{P}^1$
branched in highly singular curves.  However, 
to the authors'
knowledge, this is the first construction of such surfaces
as quasismooth hypersurfaces in a weighted projective space.

The example is defined as follows.  Consider the Delsarte polynomial
\begin{equation}
\label{eq:elliptic}
f \coloneqq x_0^{4s} + x_1^{4s} x_2 + x_2^{2s+1} x_3 + x_3^{2s+1}.
\end{equation}
This has corresponding matrix
\begin{equation*}
    A = \begin{pmatrix}
    4s & 0 & 0 & 0 \\ 
    0 & 4s & 1 & 0 \\
    0 & 0 & 2s+1 & 1 \\
    0 & 0 & 0 & 2s+1
\end{pmatrix}, \text{ with inverse }
A^{-1} = \begin{pmatrix}
    \frac{1}{4s} & 0 & 0 & 0 \\ 
    0 & \frac{1}{4s} & -\frac{1}{4s(2s+1)} & \frac{1}{4s(2s+1)^2} \\
    0 & 0 & \frac{1}{2s+1} & -\frac{1}{(2s+1)^2} \\
    0 & 0 & 0 & \frac{1}{2s+1}
\end{pmatrix}.
\end{equation*}
Using the methods from \Cref{subsect:Delsartedef}, we may compute the charge 
vector $q$ and show that $f$ defines a hypersurface
$$X = X_{4s(2s+1)^2} \subseteq \mathbb{P}((2s+1)^2,4s^2 + 2s + 1,8s^2,4s(2s+1)).$$
Here $m = 4s(2s+1)^2$ is the degree of $X$ and we'll denote the weights by
$a_0, a_1, a_2, a_3$ (in the above order).
Note that $d = m$ in this example, meaning that $X$ is a quotient of the
Fermat surface of the same degree.  The rest of this section
will be devoted to proving the properties listed in
\Cref{thm:elliptic_ex}.

\begin{claim}
For every $s \geq 1$, the hypersurface $X$ is well-formed and quasismooth.
\end{claim}

\begin{proof}
An elementary calculation confirms the following facts about the weights:
\begin{equation*}
    \gcd(a_0,a_3) = 2s+1, \gcd(a_2,a_3) = 4s, \text{ and }\gcd(a_i,a_j)
    = 1 \text{ for all other pairs $i \neq j$}.
\end{equation*}
This confirms in particular that any three distinct weights have gcd $1$, so
$\mathbb{P}$ is well-formed.  Also, $X$ does not contain any one-dimensional
stratum of $\mathbb{P}$, so $X$ is well-formed.  The equation $f = 0$ defines
a subvariety of $\mathbb{A}^4$ smooth away from the origin, meaning that
$X$ is quasismooth.
\end{proof}

Since $X$ is quasismooth, it has only cyclic quotient singularities and is 
simply connected 
\cite[Theorem 3.2.4]{Dolgachev}.

Now, we can calculate the various lattices from \Cref{subsect:quotient_setup}.
Identifying $M' \cong \Z^{4}$ and viewing $M'' \rightarrow M \rightarrow M'$ 
as sublattices, we have that $M''$ is identified with $d \Z^4$ and
the lattice $M \subseteq M'$ in between is identified with the image of the matrix

$$B^{\mathsf{T}} = d(A^{-1})^{\mathsf{T}} = \begin{pmatrix}
    (2s+1)^2 & 0 & 0 & 0 \\ 
    0 & (2s+1)^2 & 0 & 0 \\
    0 & -(2s+1) & 4s(2s+1) & 0 \\
    0 & 1 & -4s & 4s(2s+1)
\end{pmatrix}.$$

In particular, the columns of $B^{\mathsf{T}}$ generate $M \subseteq \Z^4$.
Modulo $d = 4s(2s+1)^2$, the third and fourth columns are multiples of the second,
so $M/M'' \cong \Z/4s \times \Z/(4s(2s+1)^2)$ with generators the first
and second columns $v_0$ and $v_1$,
respectively.  We next identify the subgroup $\mathcal{A} \subseteq M/M''$
(using the notation from \Cref{subsect:algorithm}).

Notice that the sum of coordinates of $v_1$ is $4s^2 + 2s + 1 = a_1$.
Consider the map $h: M/M'' \rightarrow \mathbb{Z}/d$, 
given by adding coordinates modulo $d$. 
The map $h$ sends $v_1$ to $a_1$, and since $\gcd(a_1, d) = 1$, 
it follows that $a_1$ is a generator of $\mathbb{Z}/d$. 
In particular, the map $h$ is surjective. Therefore, the group
$\mathcal{A} = h^{-1}(0)$ has index $d$ in $M/M''$, showing 
$|\mathcal{A}| = 4s$.
We claim that a generator of $\mathcal{A}$ is given by:
$$w = (2s-1)v_0 + (2s+1)^2 v_1 = \begin{pmatrix}
    (2s-1)(2s+1)^2 \\
    (2s+1)^2 \\
    (2s-1)(2s+1)^2 \\
    (2s+1)^2
\end{pmatrix}.$$
Here we have used the identities
$(2s+1)^4 \equiv (2s+1)^2 \bmod d$ and 
$-(2s+1)^3 \equiv (2s-1)(2s+1)^2 \bmod d$
to simplify, and we note that $w$ has order $4s$, hence
$\langle w \rangle = \mathcal{A}$.
This shows: $\mathcal{A} \cong \Z/4s$.
As an aside, \Cref{lem:dual_iso} proves that 
$\mathrm{Aut}_{\mathrm{tor}}(X) \cong \Z/4s$ also.
For $s \geq 2$, this means that the criterion in 
\Cref{thm:low_exp} fails.  Nevertheless, we can still
prove directly that $X_s$ has maximal Picard number.

Notice that the sum of entries of $w$ is $d$, so that 
$w \in S_0$.  The nonzero multiples of $w$ (modulo $d$) comprise
the set $\mathcal{A}^{\circ}$ for this example, 
again using the notation from 
\Cref{subsect:algorithm}.

\begin{claim}
\label{clm:maxPic}
The surface $X_s$ satisfies $p_g = s$ and has maximal Picard number.
\end{claim}

\begin{proof}
We've already observed that $w \in S_0$.  Suppose $t$ is an integer
in the set $\{1,\ldots,4s-1\}$. We claim that:
\begin{equation*}
\begin{cases}
    tw \in S_1, & \text{ if } t \text{ is even}, \\
    tw \in S_0, & \text{ if } t \text{ is odd}, t < 2s, \\
    tw \in S_2, & \text{ if } t \text{ is odd}, t > 2s.
\end{cases}
\end{equation*}
Indeed, noticing that we can divide out by $(2s+1)^2$ and that
the vector $w$
has two copies of the same pair of coordinates, the claim above
follows from the identities:
\begin{equation*}
\begin{cases}
    2m(2s-1) \bmod 4s = 4s - 2m, & \text{ if } t = 2m \text{ is even}, \\
    (2m+1)(2s-1) \bmod 4s = 2s - 2m - 1, & \text{ if } t = 2m+1 \text{ is odd}, t < 2s, \\
    (2m+1)(2s-1) \bmod 4s = 6s - 2m - 1, & \text{ if } t = 2m+1 \text{ is odd}, t > 2s.
    \end{cases}
\end{equation*}

We count that there are $s$ multiples (the odd numbers up to $2s$) lying
on $S_0$, which means that $p_g = s$.  Also,
since $4s$ is even, all the elements of $(\Z/4s)^*$ are odd. Therefore, if some
lattice point
$\alpha \in \mathcal{A}^{\circ}$ is in $S_0$ then the above 
shows that $t \alpha \in S_0$ 
or $t \alpha \in S_2$ holds for every $t \in (\Z/4s)^*$.
By \Cref{thm:Lefschetz_formula}, $X$ has maximal Picard number.
\end{proof}

When $s = 1$, this example is the quasismooth canonical K3 surface
$$X_{36} \coloneqq \{x_0^4 + x_1^4 x_2 + x_2^3 x_3 + x_3^3 = 0\} 
\subseteq \mathbb{P}(9,7,8,12)$$
with maximal Picard number $\rho = 20$.
This is a special member of family No.\ 84 on
Reid's famous list of the $95$ families of quasismooth K3 hypersurfaces
in weighted projective space (\cite[Section 4.5]{Reid79}, and see 
\cite[Section 13.3]{Iano-Fletcher} for the full list).

However, when $s \geq 2$, $X$ has worse-than-canonical singularities,
so computing plurigenera is more difficult.
To prove that $X$ is an elliptic
surface, we instead calculate the canonical map.

\begin{claim}
For each $s \geq 2$, the canonical linear system $|K_X|$ induces a 
rational map $X \dashrightarrow \mathbb{P}^{s-1}$ which is the composition
$$X \dashrightarrow \mathbb{P}^1 \rightarrow \mathbb{P}^{s-1}$$
of the pencil of monomials $x_0^2 x_2, x_1^2 x_3$ with the Veronese
embedding.  The general fiber of this pencil is a genus $1$ curve,
so in particular a resolution of $X$ is an elliptic surface
with Kodaira dimension 
$1$ for each $s \geq 2$.
\end{claim}

The same pencil of monomials also gives
the K3 surface $X$ from the $s = 1$ case
an elliptic surface structure, though of course this is no longer the
map induced by the canonical linear series.

\begin{proof}
We first note that by the adjunction formula,
$$K_X = \mathcal{O}_X(m - a_0 - a_1 - a_2 - a_3) = 
\mathcal{O}_X(16s^3 - 8s^2 - 6s - 2).$$

Factoring, $16s^3 - 8s^2 - 6s - 2 = (s-1)(16s^2 + 8s + 2)$, and we notice
that the latter factor can be written 
$16s^2 + 8s + 2 = 2a_0 + a_2 = 2a_1 + a_3$.  Therefore, the following
monomials are in the canonical linear system:
$(x_0^2 x_2)^{s-1}, (x_0^2 x_2)^{s-2}(x_1^2 x_3), \ldots,
(x_1^2 x_3)^{s-1}$.

This accounts for a dimension $s$ subset of $H^0(X,K_X)$.  Since
$p_g = s$ by \Cref{clm:maxPic}, these monomials actually span
$H^0(X,K_X)$.  This shows that the canonical map is the rational
map $\theta: X \dashrightarrow \mathbb{P}^1$ defined by 
$[x_0:x_1:x_2:x_3] \mapsto [x_0^2 x_2: x_1^2 x_3]$ composed
with the Veronese embedding $\mathbb{P}^1 \rightarrow \mathbb{P}^{s-1}$.

It remains to understand the fibers of the rational map $\theta$.  Since
$\theta$ is defined by a pencil of monomials, it is the restriction of
a toric rational map $\theta: \mathbb{P} \dashrightarrow \mathbb{P}^1$.  The
fiber $F$ of the corresponding map of tori is a two-dimensional torus.  Our
strategy will be to identify what hypersurface $X$ cuts out in this
two-dimensional torus and calculate its genus using Baker's theorem.  By focusing on 
hypersurfaces of tori, we can avoid explicitly resolving the singularities
and indeterminacies of the map $\mathbb{P} \dashrightarrow \mathbb{P}^1$.

The following figure illustrates the toric setup, where the first diagram
shows the toric maps at play, the second the corresponding cocharacter lattices,
and the last the character lattices.  Here $\mathbf{a} =
\begin{pmatrix} a_0 \\ a_1 \\ a_2 \\ a_3 \end{pmatrix}$
is the vector of weights and $\mathbf{1} = \begin{pmatrix} 1 \\ 1\end{pmatrix}$.

\begin{equation*}
\label{eq:toricsetup_elliptic}
\begin{tikzcd}[ampersand replacement=\&]
 \& \& \mathbb{A}^4 \setminus \{0\} \arrow[r,dashed] \arrow[d] \& 
 \mathbb{A}^2 \setminus \{0\} \arrow[d] \& \\
 \& F \arrow[r] \& \mathbb{P} \arrow[r, dashed, "\theta"] \& \mathbb{P}^1 \& \\
 \& \& \Z^4 \arrow[r,"{\begin{psmallmatrix} 2 & 0 & 1 & 0 \\
 0 & 2 & 0 & 1 \end{psmallmatrix}}"] \arrow[d,two heads] \& \Z^2 \arrow[d,two heads] \& \\
 0 \arrow[r] \& N_F \arrow[r] \& \Z^4/\Z  \mathbf{a} \arrow[r] 
 \& \Z^2/\Z  \mathbf{1} \arrow[r] \& 0 \\
 \& \& \Z^4 \& \arrow[l,"{\begin{psmallmatrix} 2 & 0 \\ 0 & 2 \\
 1 & 0 \\ 0 & 1 \end{psmallmatrix}}",swap] \Z^2 \& \\
 0 \& M_F \arrow[l] \& \mathbf{a}^{\perp} \arrow[l,"\eta"] \arrow[u,hook]
 \& \mathbf{1}^{\perp} \arrow[l] \arrow[u,hook] \& 0 \arrow[l] 
\end{tikzcd}
\end{equation*}

We have that $\theta$ is induced by the matrices shown because the rows (resp.
columns) correspond to the monomials $x_0^2 x_2$ and $x_1^2 x_3$.
The bottom rows of the two lattice diagrams are exact.
Now, we compute the intersection of $X$ with the torus fiber
$F = \Spec(\C[M_F])$.

The polynomial $f = x_0^{4s} + x_1^{4s} x_2 + x_2^{2s+1} x_3 + x_3^{2s+1}$ 
is a section of $\mathcal{O}_\mathbb{P}(m)$, so the points corresponding to
these four monomials sit on a translate of $\textbf{a}^{\perp} \subseteq \Z^4$.
Translating by a toric divisor in the linear system of $\mathcal{O}_\mathbb{P}(m)$
doesn't change the intersection with the torus, so we can find an explicit
Newton polytope for $X$ inside the lattice $\textbf{a}^{\perp}$
corresponding to the torus of $\mathbb{P}$ by translating
the above four points to lie in $\textbf{a}^{\perp}$.  Explicitly, we'll translate
by $x_0^{-4s}$ to get the polytope in $\textbf{a}^{\perp}$ generated by the 
four points
$$p_1 = \begin{pmatrix} 0 \\ 0 \\ 0 \\ 0 \end{pmatrix}, 
p_2 = \begin{pmatrix} -4s \\ 4s \\ 1 \\ 0 \end{pmatrix},
p_3 = \begin{pmatrix} -4s \\ 0 \\ 2s+1 \\ 1 \end{pmatrix},
p_4 = \begin{pmatrix} -4s \\ 0 \\ 0 \\ 2s+1 \end{pmatrix}.$$

Finally, the Newton polytope of the equation of the hypersurface $X \cap F 
\subseteq F$ is given by mapping the points $p_1,\ldots,p_4$ to $M_F$ via
the surjection $\eta: \mathbf{a}^{\perp} \rightarrow M_F$ defined to be
the cokernel
of $\mathbf{1}^{\perp} \rightarrow \mathbf{a}^{\perp}$.
Define $\bar{p}_i \coloneqq \eta(p_i)$ for each $i = 1,\ldots,4$.
The following lemma gives the property we need of the image polytope.

\begin{lemma}
\label{lem:onelatpoint}
For each $s \geq 1$, the lattice quadrilateral $Q$ defined by the convex hull of 
$\bar{p}_1,\bar{p}_2,\bar{p}_3,\bar{p}_4$ in $M_F$ has exactly 
one interior lattice point.
\end{lemma}

\begin{proof}
The image of $\textbf{1}^{\perp} = \Z  \begin{pmatrix} 1 \\ -1 \end{pmatrix}
$ under the map
$\textbf{1}^{\perp} \rightarrow \textbf{a}^{\perp}$ is the rank $1$ sublattice
generated by 
$$u = 1\begin{pmatrix} 2 \\ 0 \\ 1 \\ 0 \end{pmatrix} 
- 1\begin{pmatrix} 0 \\ 2 \\ 0 \\ 1 \end{pmatrix}
= \begin{pmatrix} 2 \\ -2 \\ 1 \\ -1 \end{pmatrix}.$$
We can view the quotient $M_F = \textbf{a}^{\perp}/\Z  u$ as sitting inside
$\Z^4/\Z  u$, 
the latter of which we identify with $\Z^3$ by choosing the basis $\bar{e}_1,\bar{e}_2,
\bar{e}_3$ corresponding to the images of the standard basis vectors in $\Z^4$.
We can do this since $e_1,e_2,e_3,u$ form a basis of the original $\Z^4$ (the matrix with
these four columns has determinant $-1$).

In these new coordinates, we calculate images of the $p_i$ by writing them in terms of the
above basis to get

$$\bar{p}_1 = \begin{pmatrix} 0 \\ 0 \\ 0\end{pmatrix}, 
\bar{p}_2 = \begin{pmatrix} -4s \\ 4s \\ 1 \end{pmatrix},
\bar{p}_3 = \begin{pmatrix} -4s+2 \\ -2 \\ 2s+2 \end{pmatrix},
\bar{p}_4 = \begin{pmatrix} 2 \\ -4s-2 \\ 2s+1 \end{pmatrix}.$$

We just need to show that the quadrilateral $Q$ determined by these points in the plane
they span has exactly one interior lattice point.  Indeed, observe that 
$\bar{p}_3 = \bar{p}_2 + \bar{p}_4$, and that $\bar{p}_3/2$ is an interior
lattice point of $Q$.  We will argue it is the only one.  We compute that the
$2 \times 2$ minors of the matrix
$$\begin{pmatrix} \bar{p}_2 & \frac{\bar{p}_3}{2}\end{pmatrix} = 
\begin{pmatrix} -4s & -2s+1 \\ 4s & -1 \\ 1 & s+1 \end{pmatrix}$$
are $8s^2, 4s^2 + 4s + 1$, and $-4s^2 - 2s - 1$.  Since 
$\gcd(8s^2,4s^2 + 4s + 1,-4s^2 - 2s -1) = 1$, these two vectors
extend to some basis of $\Z^3$ \cite[Lemma 2, p. 15]{Cassels} and 
in particular generate the sublattice 
in the plane they span.  Hence after one last coordinate change, the 
four corners $\bar{p}_1,\ldots,\bar{p}_4$ correspond in the basis
$\bar{p}_2, \bar{p}_3/2$ to
the points $\begin{pmatrix} 0 \\ 0\end{pmatrix}, 
\begin{pmatrix} 1 \\ 0 \end{pmatrix}, 
\begin{pmatrix} 0 \\ 2 \end{pmatrix}, 
\begin{pmatrix} -1 \\ 2\end{pmatrix}$ in $\Z^2$, respectively.  These
indeed define a quadrilateral $Q$ enclosing exactly $1$ interior 
point, namely $\begin{pmatrix} 0 \\ 1 \end{pmatrix}$.\end{proof}

Since we can take the coefficients of the equation $f$ defining $X$ to be
general, the hypersurface $X \cap F$ is Newton non-degenerate.
Therefore, Baker's theorem and \Cref{lem:onelatpoint}
show that the unique smooth,
projective curve
isomorphic to the fiber $X \cap F$ of the rational map 
$\theta: X \dashrightarrow \mathbb{P}^1$ has genus $1$.

The canonical pencil is a birational invariant of klt varieties, so after
suitably blowing up $X$ to resolve singularities and the pencil $\theta$,
we get a morphism $\widetilde{X} \rightarrow \mathbb{P}^1$ whose general
fiber is a smooth genus $1$ curve.  This proves that $\widetilde{X}$ is an 
elliptic surface.  Since $p_g = s \geq 2$, it has Kodaira dimension $1$.
\end{proof}

\begin{remark}
Once the Lefschetz number of an elliptic surface $X$ 
over $\mathbb{P}^1$
is known, it is possible
to algorithmically compute the rank of the Mordell-Weil group
of $X$ using the \textit{Shioda-Tate formula} \cite{Shioda72,Tate}.
This has been implemented by Heijne \cite{Heijne12} to 
systematically compute the Mordell-Weil ranks for
elliptic curves over $\C(t)$
``of Delsarte type".  In particular,
\cite[Section 6]{Heijne12} computes the maximal 
Mordell-Weil rank
for elliptic curves with a particular Newton polygon.
By the proof of \Cref{lem:onelatpoint}, the Newton polygon
for the fiber of $X_s \dashrightarrow \mathbb{P}^1$ is isomorphic
to the quadrilateral with vertices $(0,0), (1,0), (0,2)$, and
$(-1,2)$.  After a coordinate change, this matches
the 10th polygon in Heijne's table.  As the table shows,
the Mordell-Weil rank equals $0$ for any Delsarte elliptic
curve over $\C(t)$ with this Newton polygon.
This proves that
all the examples $X_s$, despite having maximal Picard
rank, have the smallest possible Mordell-Weil rank.
\end{remark}

A computer search finds many similar infinite
sequences of examples with maximal Picard number
and increasing
$p_g$, but all of those we found
are sequences of elliptic surfaces with similar properties
to those defined by \eqref{eq:elliptic}.
In contrast, general type examples like the one
in \Cref{ex:gen_type_low_exp} seem to be sporadic
and very rare.

\end{document}